\documentclass[12pt,a4paper]{article}

\usepackage{color} 
\usepackage{amsmath}
\usepackage{amsfonts}
\usepackage{amssymb}
\usepackage{amsthm} 
\usepackage{graphicx} 
\graphicspath{{./figures/}} 

\usepackage{mathtools}

\newcommand{\Spin}{{\mathrm{Spin}}}
\newcommand{\SL}{{\mathrm{SL}}}
\newcommand{\PSL}{{\mathrm{PSL}}}
\newcommand{\GL}{{\mathrm{GL}}}

\newcommand{\U}{{\mathrm U}}
\newcommand{\SO}{{\mathrm{SO}}}

\newcommand{\dd}{{\mathrm d}}  

\newcommand{\R}{{\mathbb R}}   
\newcommand{\C}{{\mathbb C}}   
\newcommand{\Z}{{\mathbb Z}}   

\DeclareMathOperator\End{End}
\DeclareMathOperator\Hom{Hom}
\DeclareMathOperator\tr{Tr}

\newcommand\pd[2]{{\frac{\partial #1}{\partial #2} }} 

\newcommand{\git}  
{\mathbin{   
  \mathchoice{/\mkern-6mu/}
    {/\mkern-6mu/}
    {/\mkern-5mu/}
    {/\mkern-5mu/}}}

\theoremstyle{definition}
\newtheorem{lemma}{Lemma} 
\newtheorem{definition}{Definition} 
\newtheorem{example}{Example} 
 
\newtheorem*{example*}{Example} 

 \newcommand{\q}{{\mathfrak q}}   
 \newcommand{\Q}{q}  
 \newcommand{\X}{Q} 
\newcommand{\e}{{\mathfrak e}} 

\newcommand{\lBracket}{[\![}  
\newcommand{\rBracket}{]\!]}  

\begin{document}

\title{Finite spectral triples for the fuzzy torus}

\author{John W. Barrett
\thanks{Corresponding author} 
\\ \\
School of Mathematical Sciences\\
University of Nottingham\\
University Park\\
Nottingham NG7 2RD, UK\\
\and James Gaunt\thanks{Current address:  Department of Mathematics, Heriot-Watt University, Colin Maclaurin Building, Riccarton, Edinburgh EH14 4AS, UK}
\\ \\
School of Mathematical Sciences\\
University of Nottingham\\
University Park\\
Nottingham NG7 2RD, UK\\
\\
\\
E-mail john.barrett1@nottingham.ac.uk
\\
j.gaunt@hw.ac.uk}

\date{30 October 2024}

\maketitle

\begin{abstract} Finite real spectral triples are defined to characterise the non-commutative geometry of a fuzzy torus. The geometries are the non-commutative analogues of flat tori with moduli determined by integer parameters. Each of these geometries has four different Dirac operators, corresponding to the four spin structures on a torus. The spectrum of the Dirac operator is calculated. It is given by replacing integers with their quantum integer analogues in the spectrum of the corresponding commutative torus. 
\end{abstract}

\section{Introduction}
A Riemannian spin geometry can be expressed in algebraic terms using the algebra of smooth functions on the manifold and the Dirac operator on the spinor bundle.  This point of view allows for a significant generalisation of geometry by removing the restriction that the algebra is commutative, adjusting the axioms in as natural a way as possible. The mathematical structure that results is called a real spectral triple \cite{Connes:1995tu} and encompasses both commutative and non-commutative geometries.

Among the non-commutative geometries, it is of significant interest to examine examples that are analogues of Riemannian spin manifolds, or, more specifically, approximations of Riemannian spin manifolds. This paper contributes to this by defining and studying the properties of some non-commutative analogues of the flat torus called fuzzy tori. These analogues are all \emph{finite} in the sense that the algebra is a (finite-dimensional) matrix algebra, and the Dirac operator is also a matrix. This notion of finiteness turns out to be richer than the corresponding notion in commutative geometry, which would just be the study of finite sets of points. 

The flat torus is obtained by identifying the opposite sides of a parallelogram in Euclidean space. While this is simple as a geometry, it provides interesting examples due to the free parameters in the shape of the torus and the four possible spin structures. Both the geometry and the spin structure can be read off from the spectrum of the Dirac operator. 

A fuzzy torus is a finite non-commutative analogue of a flat torus and its geometry is specified by a choice of two elements of its algebra.
The main result of this paper is the construction of a Dirac operator for this fuzzy torus geometry, forming a real spectral triple. The general formula is analogous to the commutative formula for a Dirac operator constructed with a particular frame field. In simple examples, where the fuzzy torus algebra is a simple matrix algebra, the spectrum of the Dirac operator is a deformation of the spectrum of a commutative Dirac operator and one can determine the corresponding spin structure of the commutative torus from the spectrum. This coincides with the spin structure determined by the frame field.

In further examples with a non-simple algebra, it is possible to construct Dirac operators that correspond to all four possible spin structures on the commutative torus. This is done by constructing the non-commutative analogue of a four-fold covering of the torus.

\paragraph{The torus} The algebra of functions on the torus $T^2$ is the commutative algebra generated by 
\begin{equation} U=e^{i\theta},\quad V=e^{i\phi}.
\end{equation}
Defining 
\begin{equation} X=e^{i(a\theta+b\phi)},\quad Y=e^{i(c\theta+d\phi)},
\end{equation}
for integers $a,b,c,d$, the principal objects of study are differential operators specified using $X$ and $Y$. These correspond to a Riemannian metric on the torus. In the simplest case $X=U$ and $Y=V$, this is the square metric
\begin{equation}\dd\theta^2+\dd\phi^2.
\end{equation}
In the general case, the metric is the pull-back of this by the covering map
\begin{equation}(\theta,\phi)\mapsto (a\theta+b\phi,c\theta+d\phi). 
\end{equation}

The finite non-commutative torus is the study of unitary matrices $U$ and $V$, subject to the relation
\begin{equation} UV=\Q VU
\end{equation}
for a complex number $\Q$, which has to be a root of unity. This is simultaneously a deformation of the torus and a truncation to finite dimensions. 

The \emph{fuzzy torus} is the study of differential operators on this algebra that are constructed using monomials $X$ and $Y$ in the preferred `coordinates' $U$ and $V$. These operators are the non-commutative analogues of the corresponding torus metrics. 

The simplest example of the relation is given by $U=C$ and $V=S$, the well-known clock and shift matrices (see Example \ref{ex:cs}),  which generate the algebra $M_N(\C)$ of all $N\times N$ matrices. More generally, one can take $U$ and $V$ to be monomials in $C$ and $S$. 
A particularly important example is given by $U=C^2$ and $V=S^2$. The algebra generated by $U$ and $V$ is a subalgebra of $M_N(\C)$, giving a non-commutative version of the regular covering map $(\theta,\phi)\mapsto (2\theta,2\phi)$. The characters of the deck transformation group $\Z_2\times\Z_2$ distinguish the different spin structures on the base of the covering. In the non-commutative case, the analogous construction allows one to define fuzzy tori with different spin structures.

\paragraph{Spectral triples} The idea of a spectral triple is that it provides a mathematical framework for a Dirac operator \cite{ConnesBook}.
A spectral triple consists of an algebra $\mathcal A$, a Hilbert space $\mathcal H$ on which $\mathcal A$ acts, and a Dirac operator $D\colon \mathcal H\to \mathcal H$. For a Riemannian $\text{spin}$ manifold $\mathcal{M}$, $\mathcal A$ is the algebra of smooth functions on $\mathcal{M}$,  $\mathcal H$ is the square-integrable sections of the $\text{spin}$ bundle, and $D$ is the usual Dirac operator. The spin structure defines an antilinear map $\mathcal J\colon\mathcal H\to \mathcal H$, called the real structure. 

In the non-commutative generalisation, $\mathcal A$ is obviously a non-com\-mu\-ta\-tive algebra. There is a distinct right action of $ \mathcal{A}$ on $\mathcal H$, promoting it to a bimodule rather than just a left module. The real structure plays the role of interchanging the left and right actions.

\paragraph{The results} The main result of this paper is the construction of finite real spectral triples that are non-commutative analogues of various metrics for the torus, for all four spin structures. This is done by extending the construction of the fuzzy torus. 

The Dirac operator is constructed using algebra elements $X$ and $Y$ in Definition \ref{def:spectral}. It is a very close non-commutative analogue of the commutative one, with derivatives replaced by commutators and functions replaced by anticommutators.  The non-trivial part of the construction, in the commutative case, is that it is necessary to use a \emph{rotating frame} in the description of the Dirac operator. In this case, the spin connection coefficients are non-zero. 
This is necessary because there is no matrix analogue of the vector fields $\partial/\partial\theta$ and $\partial/\partial\phi$ that would have to be used for a non-rotating frame \cite{LLS}. For example, a matrix analogue of the vector field $\partial/\partial\theta$ 
for the algebra generated by clock and shift matrices $C$ and $S$
would be a matrix $L$ satisfying 
\begin{equation} [L,C]=iC,\quad\quad [L,S]=0.\end{equation}
However, no such matrix exists.

The non-commutative version of the Dirac operator has a term that is a direct analogue of the connection coefficients. The role of the rotating frame is somewhat mysterious from the conceptual point of view, but it is worth noting that it is one of the factors that determines the spin structure. 

The eigenvalues of the Dirac operator are calculated and it is found, rather beautifully, that these are the $q$-number analogues of the eigenvalues for the commutative case, a result already known \cite{SS} for the rather simpler Laplace operator \cite{Kimura:2001uk}. In fact, the spectrum of the fuzzy Dirac operator is exactly the set of square roots of the spectrum of the corresponding fuzzy Laplacian. It is worth noting that the correct definition of the non-commutative analogues of the connection coefficients is essential to obtain this result. Note that the construction of a Dirac operator on a fuzzy torus in \cite{SS} was limited to the derivative (commutator) part only and does not appear to give a spectrum that corresponds to the commutative case.

The results are quite different to existing constructions of the Dirac operator on the rational non-commutative torus.  The algebra of the rational non-commutative torus can be understood as the space of sections of a bundle of matrix algebras over the torus, as explained in \cite{H-KS} and \cite[Section 12.2]{G-BVF}. 
Dirac operators on the rational non-commutative torus are constructed in \cite{Paschke:2006ooy}; see also \cite{Venselaar:2010ldp} for further details. They each have the same spectrum as a Dirac operator in the commutative case and the linear growth of eigenvalues is rather different to the periodic $q$-number spectra of the finite spectral triples presented here. The different spin structures on the rational non-commutative torus are understood in terms of two-fold coverings in \cite{Carotenuto:2018cmx}, using a different covering for each spin structure rather than the universal four-fold covering developed here. However, the overall idea of using a non-trivial character on a covering space to change spin structure is the same.

\paragraph{Summary} The relevant geometry of a flat torus is summarised in Section \ref{sec:torus}. The non-commutative torus is described in Section \ref{sec:nctorus}, restricted to the finite case, including the analogues of translations of the coordinates and the modular transformations. There are also some heuristics for the commutative limit that guide the ideas in later sections.

The fuzzy torus is introduced in Section \ref{sec:fuzzy}. The notion of a fuzzy space is slightly more sophisticated than that usually considered in the literature. Instead of considering just an algebra with distinguished generators, the fuzzy torus also has a bimodule over the algebra. Of course this includes the standard case where the bimodule is the algebra itself, but it also allows for a significant generalisation where the bimodule is different, corresponding to non-trivial line bundles over the torus with $\Z_2$ twists. These non-trivial bundles are the precursors of different spin structures in later sections.
A definition of the Laplace operator is given and its properties are studied. In particular, it is shown that with suitable choices of matrices $X$ and $Y$, the Laplace operator is the analogue of the commutative Laplace operator on a torus with a flat metric. The geometry of the torus can be seen from studying the spectrum of the fuzzy Laplacian.

The Dirac operator on the flat torus is studied in Section \ref{sec:dirac} and specialised to the torus in Section \ref{sec:diractorus}. These sections summarise the formalism of the Dirac operator in a very explicit way. While the usual abstract notation in differential geometry is very concise and efficient, it hides the fact that there are many equivalent explicit formulas. Since the symmetries of the commutative and non-commutative cases are different, it must be that not all explicit formulas in the commutative case have non-commutative analogues. Thus it is necessary to present the formulas in a very particular way to allow the correct generalisation to the non-commutative case. 

The main results are given in Section \ref{sec:fuzzydirac}. The finite real spectral triple is defined in Definition \ref{def:spectral}. This is first examined for the square torus, with one particular spin structure, and then generalised for other geometries of flat tori, illustrated by some plots of the spectra. Finally, it is then shown how to generalise the results to all the four spin structures on a given fuzzy torus using the non-commutative analogue of a four-fold covering map. It is explained how this construction appears to be a non-commutative and discrete analogue of Marsden-Weinstein symplectic reduction.

\paragraph {Acknowledgements}  Thanks are due to the hospitality of the Mathematisches Forschungsinstitute Oberwolfach during July 2018, where this work was presented, as summarised in the meeting report \cite{MFO}. Part of this work also appears in JG's PhD thesis \cite{PhD} and the support of the University of Nottingham is acknowledged.
\smallskip

Declarations of interest: none

\section{The flat torus}\label{sec:torus}
The torus is defined as the quotient
 \begin{equation} T^2=\R^2/(2\pi\Z)^2.
\end{equation}

Let $A\colon\R^2\to\R^2$ be a linear map with integer coefficients, i.e., $A\in M_2(\Z)$. It maps lattice points to lattice points and hence determines a map $\widetilde A\colon T^2\to T^2$. 

If $A$ is invertible and $A^{-1}$ also has integer coefficients, then $\widetilde A$ is invertible, and the set of all such $A$ is written $\GL(2,\Z)$.  Note that in this case, $\det A$ and its inverse are integers, so $\det A=\pm1$. 
If in addition $\det A=1$, the transformation is orientation-preserving and is called a modular transformation. The set of all elements with determinant $1$ defines the group $\SL(2,\Z)$, which will be called the modular group (note that this name is also commonly used for the quotient $\PSL(2,\Z)=\SL(2,\Z)/\{\pm1\}$).

A flat metric on $\R^2$ is determined by the formula
\begin{equation}g=g_{\mu\nu}\,\dd x^\mu\dd x^\nu
\end{equation}
 using the Einstein summation convention
with constant coefficients $g_{\mu\nu}$ and coordinates $(x^1,x^2)$. All such metrics that are positive-definite are obtained by the pull-back of the standard Euclidean metric $\delta=(\dd x^1)^2+(\dd x^2)^2$ by a linear map $A\in\GL(2,\R)$,
\begin{equation} g=A^\bullet\delta\label{eq:pbmetric}\end{equation}
and the coefficients are given by
\begin{equation} g_{\mu\nu}=A_\mu^\rho A_\nu^\sigma\delta_{\rho\sigma}.\end{equation}
Note that the pull-back is denoted $A^\bullet$ rather than the usual $A^*$ so as not to confuse with the adjoint of the matrix $A$. The corresponding metric on the torus will be denoted 
$\widetilde g$ and is called a \emph{flat torus}. A flat torus has zero Riemannian curvature but, moreover, the coordinates are chosen so that the metric coefficients are constant.

The following lemma describes when $\widetilde g$ is equivalent to a second flat metric, $\widetilde g'$.
\begin{lemma} The flat tori $(T^2,\widetilde g)$ and $(T^2,\widetilde g')$ are isometric iff there exists $B\in\GL(2,\Z)$ such that $g'=B^\bullet g$.
\end{lemma}
\begin{proof} Suppose $B\in\GL(2,\Z)$ is such that $g'=B^\bullet g$. Since $B\Z^2\subset\Z^2$, it determines a mapping on the quotient, $\widetilde B\colon T^2\to T^2$, and $\widetilde g'=\widetilde B^\bullet \widetilde g$.

Conversely, suppose the two flat tori are isometric by $\widetilde B\colon T^2\to T^2$. Then there is an  isometry $B\colon\R^2\to \R^2$, taking $g$ to $g'$, that projects to $\widetilde B$. This isometry is an affine map, and by composing by a translation one can assume it is linear, i.e., $B0=0$. Then $B$ maps lattice points to lattice points, and so $B\in M_2(\Z)$. Similarly, $B^{-1}\in M_2(\Z)$, so $B\in\GL(2,\Z)$.
\end{proof}

The standard metric on the torus, $\widetilde\delta$, is obtained by identifying the edges of a square with side length $2\pi$, and so is called a square torus. If $A$ has integer coefficients, $A\in\GL(2,\R)\cap M_2(\Z)$, then the metric is obtained by the pull-back of $\widetilde\delta$ by the map $\widetilde A$. These metrics will play an important role in the fuzzy torus and will be called \emph{integral}.

If $g=A^\bullet\delta$ then $B^\bullet g=(AB)^\bullet\delta$. Thus the isometry class of integral metrics on the torus is determined by the equivalence class of matrices $\{AB|B\in\GL(2,\Z)\}$. Each equivalence class has a unique representative called the Hermite normal form of the integer matrix $A$. These are the matrices
\begin{equation}A=\begin{pmatrix}a&0\\c&d\end{pmatrix}
\end{equation} 
with $a>0$ and $0\le c<d$. These matrices will be useful for constructing examples.

\subsection{Poisson structure}
The orientation of the torus can be registered with a 2-form that is compatible with the metric, making a K\"ahler structure. For the standard Euclidean metric, define
\begin{equation} \omega_0=\;\dd x^1\wedge\dd x^2.
\end{equation}
and the corresponding 2-form on the torus is denoted $\widetilde\omega_0$.
The standard orientation $\mathcal O$ is the orientation of $T^2$ such that
\begin{equation}
\int_{T^2, \mathcal O}\widetilde\omega_0=4\pi^2,
\end{equation}
i.e., the integral is positive.

For a general metric $g$, given by \eqref{eq:pbmetric}, the compatible 2-form on $\R^2$ is the pull-back $\omega= A^\bullet\omega_0$, and determines a 2-form $\widetilde\omega$ on the torus. Integrating this gives
\begin{equation}\int_{T^2, \mathcal O} \widetilde\omega=(\det A)\, 4\pi^2=\pm\int_{T^2}\sqrt{\det g}\;\dd x^1\dd x^2.
\end{equation}
The sign in this equation can be taken to indicate an orientation of the torus, positive being the standard orientation $\mathcal O$ and negative the opposite one, $-\mathcal O$. Note that for an integral metric the factor $(\det A)$ counts the winding number of the map $\widetilde A$.

Another approach to specifying the metric is to give coordinates in which the metric components take a standard form. It is convenient to use the explicit notation
\begin{equation} x^1=\theta,\quad x^2=\phi.\end{equation}
 Since these are not functions on the torus, one takes exponentials  instead. The standard metric and orientation is specified by the ordered pair of functions
\begin{equation} U=e^{i\theta},\quad V=e^{i\phi}.
\end{equation}
 so that 
\begin{equation}\widetilde\delta=\dd\overline U\dd U+\dd\overline V\dd V,\quad\quad \widetilde\omega_0=-\overline U\overline V\,\dd U\wedge\dd V.\label{eq:cood-defs}
\end{equation}
For a general integral metric, if $A\in \GL(2,\R)\cap M_2(\Z)$ is the matrix
\begin{equation}A=\begin{pmatrix}a&b\\ c&d\end{pmatrix}\label{eq:A}
\end{equation}
 the pull-back functions are
\begin{equation} X=e^{i(a\theta+b\phi)}=U^aV^b, \quad Y=e^{i(c\theta+d\phi)}=U^cV^d \label{eq:transform}
\end{equation}
and the pull-backs of \eqref{eq:cood-defs} now determine $\widetilde g$ and $\widetilde\omega$,
\begin{gather}\widetilde g=\dd\overline X\dd X+\dd\overline Y\dd Y=(a^2+c^2)\dd\theta^2+2(ab+cd)\dd\theta\dd\phi+(b^2+d^2)\dd\phi^2\label{eq:cood-defs-gen}
\\
\widetilde \omega=-\overline X\overline Y\,\dd X\wedge\dd Y=(ad-bc)\dd\theta\wedge\dd\phi.
\end{gather}

The  Poisson bracket is the inverse of the 2-form $\omega$. The 2-form is expressed as a tensor, using the (not universal) convention
\begin{equation} \omega=\frac12\,\omega_{\mu\nu}\,\dd x^\mu\wedge\dd x^\nu=\omega_{\mu\nu}\,\dd x^\mu\otimes\dd x^\nu,
\end{equation}
assuming $\omega_{\mu\nu}=-\omega_{\mu\nu}$. Then defining the coefficients $\Omega^{\mu\nu}$ to be the  inverse matrix, so that $\Omega^{\mu\nu}\omega_{\nu\rho}=\delta^\mu_\rho$, the Poisson bracket is
\begin{equation}\Omega(f,g)=\Omega^{\mu\nu}\,\pd f{x^\mu}\pd g{x^\nu}.
\end{equation}
For the Poisson bracket determined by $\omega_0$ this formula is
\begin{equation}\Omega_0(f,g)=\left(\pd f\phi\pd g\theta-\pd f\theta\pd g\phi\right).
\end{equation}
Applying this to exponentials gives
\begin{equation}\Omega_0(U^mV^n,U^kV^l)=(ml-nk)U^{m+k}V^{n+l}\label{eq:poissonexp}
\end{equation}
In general, one has
\begin{equation}\Omega(U^mV^n,U^kV^l)=\frac{ml-nk}{ad-bc}U^{m+k}V^{n+l}
\end{equation}
and in particular,
\begin{equation}\Omega(X,Y)=X Y.
\end{equation}

\section{The non-commutative torus}\label{sec:nctorus}

\begin{definition}
A non-commutative torus is a pair of unitary operators $U$ and $V$ on a Hilbert space $\mathfrak h$ and a complex number $\Q$ satisfying
\begin{equation} UV=\Q VU.
 \label{eq:uv}
\end{equation}
A  finite non-commutative torus is one for which $\mathfrak h$ is finite-dimensional.
\end{definition}

The unitary condition implies (by taking the adjoint of both sides) that $|\Q|=1$.  Heuristically, one can view \eqref{eq:uv} as determining a deformation of the commutative algebra of functions on a torus.

Two such tori are regarded as equivalent if there is a unitary isomorphism $P\colon\mathfrak h\to\mathfrak h'$ transforming one to the other,
\begin{equation} PUP^{-1}=U',\quad  PVP^{-1}=V'.\label{eq:trans}
\end{equation}
A torus is irreducible if there is no invariant proper subspace of $\mathfrak h$.

This paper considers exclusively the case of finite non-commutative tori, i.e., $U$ and $V$ are matrices. The basic properties were derived by Weyl \cite[IV D]{Weyl}, and Mackey \cite{Mackey} in more generality.  

If $\psi$ is an eigenvector of $U$ with eigenvalue $\lambda$, then $V\psi$ has eigenvalue $\Q \lambda$. Since there are only a finite number of eigenvalues, it must be that 
\begin{equation} 
\Q^N=1 
\end{equation}
 for some positive integer $N$. The order of  $\Q$  is the smallest possible value of $N$ and then $\Q$ is said to be a primitive $N$-th root of unity.

The Hilbert space decomposes into a direct sum of irreducibles. This is because the orthogonal complement of an invariant subspace is also an invariant subspace. Therefore to understand the general case it suffices to look at the irreducibles. 

So now assume $\mathfrak h$ is irreducible. The operators $U^N$ and $V^N$ commute with each other, and with $U$ and $V$, so in the $\lambda$ eigenspace of $U$ there are simultaneous eigenvectors of both $U^N$ and $V^N$. Let $\psi$ be one of these.
Then the subspace with basis
\begin{equation}\label{eq:basis}\{\psi, V\psi, V^2\psi,\ldots,V^{N-1}\psi\}
\end{equation}
 is an invariant subspace for $U$ and $V$. So it follows that in an irreducible the eigenspaces of $U$ are all one-dimensional and
\begin{equation}U^N=\alpha 1,\quad V^N=\beta 1,\quad \alpha,\beta\in \C.
\end{equation}
Unitarity implies $|\alpha|=|\beta|=1$. The numbers $\Q, \alpha, \beta$ characterise the irreducible geometry.

\begin{example} \label{ex:cs} The standard example is given by the  clock and shift matrices of size $N$ determined by a complex number $\q$ of order $N$, so that $\q^N=1$. These are the $N\times N$ matrices
\begin{equation}C=\begin{pmatrix}1&0&\ldots&0\\0&\q&\ldots&0\\0&0&\ddots&0\\0&0&0&\q^{N-1}\end{pmatrix},
\quad S=\begin{pmatrix}0&0&\ldots&1\\1&0&\ldots&0\\0&\ddots&0&0\\0&0&1&0\end{pmatrix}
\end{equation}
These satisfy the relations
\begin{equation}CS=\q SC,\quad C^N=1,\quad S^N=1\label{eq:CSrelations}
\end{equation}
So putting $U=C$, $V=S$, $\Q=\q$, one finds that $\alpha=\beta=1$.
\end{example}

The general irreducible case is obtained from the clock and shift matrices by a rescaling.
\begin{lemma}[Weyl] \label{lem:basis} Let $U$, $V$ be an irreducible finite non-commutative torus. Define the clock and shift matrices $C$ and $S$ according to Example \ref{ex:cs}, with $\q=\Q$. Then there is a basis of $\mathfrak h$ such that $U=\alpha^{1/N}C$, $V=\beta^{1/N}S$.
\end{lemma}

\begin{proof}
Define $C=\alpha^{-1/N}U$, $S=\beta^{-1/N}V$, for an arbitrary choice of the $N$-th roots. 
Using the following basis
\begin{equation}\label{eq:basis2}\{\psi, S\psi, S^2\psi,\ldots,S^{N-1}\psi\},
\end{equation}
the matrices are those given in Example \ref{ex:cs}.
\end{proof}

The linear maps of $\mathfrak h$ to itself,  $\End(\mathfrak h)$, form a $*$-algebra over $\C$. The $*$ operation is the adjoint of a linear map.
\begin{definition}
Let $(U,V,\mathfrak h)$ be a finite non-commutative torus. 
The torus algebra  $\langle U,V\rangle$ is the $*$-subalgebra of $\End(\mathfrak h)$ generated by $U$ and $V$. 
\end{definition}

The torus algebra is semisimple and hence isomorphic to a direct sum of matrix algebras. Let $N$ be the order of $q$. In the irreducible case, Lemma \ref{lem:basis} implies that $\langle U,V\rangle\cong M_N(\C)$. One can also see directly that the matrices $\{C^nS^m|m,n=0\ldots N-1\}$ are linearly independent and so form a basis of $M_N(\C)$. In the general case, $\langle U,V\rangle$ is isomorphic to the direct sum of a finite number of copies of  $M_N(\C)$.

In the following, the matrices $C$ and $S$ are always the clock and shift matrices defined in Example \ref{ex:cs}, whereas $U$ and $V$ are operators defined in a variety of different ways. There are interesting examples obtained by taking monomials in $C$ and $S$.

\begin{example}\label{ex:Csquard} Let $C$, $S$ be the clock and shift matrices of size $N$ with parameter $\q$. Define $U=C^2$, $V=S$. Then $UV=\Q VU$ with $\Q=\q^2$. If $N$ is odd, then $\Q$ is again order $N$ and the geometry is irreducible.

However if $N$ is even, then the order of $\Q$ is $N/2$. The geometry determined by $U$ and $V$ reduces into the direct sum of two irreducible geometries according to the eigenvalues $\beta=\pm1$ of $S^{N/2}$. The algebra $\langle U,V\rangle$ is isomorphic to $M_{N/2}(\C)\oplus M_{N/2}(\C)$.
\end{example}

It is worth comparing the finite case with an infinite non-commutative torus with the same value of $\Q$. This is called a rational non-commutative torus. It is shown in \cite[Prop 12.2]{G-BVF} that the universal torus algebra is isomorphic to a bundle of matrix algebras over a commutative torus. This commutative torus is generated by $U^N$ and $V^N$ and each fibre of the bundle is isomorphic to $M_N(\C)$. Any other torus algebra (finite or infinite) is a quotient of this universal torus algebra. In the finite case, the commutative torus is replaced by a finite subset of points determined by the eigenvalues of $U^N$ and $V^N$, so that the algebra is a direct sum of matrix algebras as described above.

\subsection{Normalised monomials}
For any torus algebra  $\langle U,V\rangle$ and a choice of square root $\Q^{1/2}$ one can define a convenient normalisation of the  monomials of the generators that  resolves the ambiguity between $V^nU^m$  and $U^mV^n={\Q}^{mn}V^nU^m$.

\begin{definition} The \emph{normalised monomials} $E^{(m,n)}$ are defined for integers $m,n$ by
\begin{equation}E^{(m,n)}=\Q^{-mn/2}U^mV^n=\Q^{mn/2}V^nU^m.
\end{equation}
\end{definition}

 The multiplication of the monomials is
\begin{equation} E^{(m,n)}E^{(k,l)}=\Q^{(ml-nk)/2} E^{(m+k,n+l)},\label{eq:Emult}
\end{equation}
so in particular,
\begin{equation}(E^{(m,n)})^k=E^{(km,kn)}
\end{equation} 
for any integer $k$.

A similar discussion of these operators can be found in \cite{D'Andrea} and are often referred to as the Weyl operators within the literature. The key distinction between the operators introduced here and those within \cite{D'Andrea} is the normalisation factor which plays a key role in the calculation of the spectrum of the scalar Laplace and Dirac operators on the fuzzy torus.

The anticommutator of two matrices is $\{A,B\}=AB+BA$. For the normalised monomials this is
\begin{equation} \{E^{(m,n)},E^{(k,l)}\}= (\Q^{(ml-nk)/2}+\Q^{-(ml-nk)/2}  ) E^{(m+k,n+l)}.
\end{equation}
which gives a commutative but non-associative product.

The commutators are given by the following expression.
\begin{equation}\begin{aligned}\frac1{\Q^{1/2}-\Q^{-1/2}} [E^{(m,n)},E^{(k,l)}]&=\frac{ \Q^{(ml-nk)/2}-\Q^{-(ml-nk)/2}}{\Q^{1/2}-\Q^{-1/2}}   E^{(m+k,n+l)}.\\
&=[ml-nk]_\Q\;E^{(m+k,n+l)},
\end{aligned}\label{eq:fuzzycommutators}\end{equation}
using the \emph{quantum integer}
\begin{equation}[n]_\Q=\frac{ \Q^{n/2}-\Q^{-n/2}}{\Q^{1/2}-\Q^{-1/2}}\label{eq:quantuminteger}
\end{equation}
defined for any integer $n$. Suppose $\Q=e^{2\pi i K/N}$ with $K$ coprime to $N$, then $\Q^{1/2}=\pm e^{\pi i K/N}$ and
\begin{equation}
[n]_\Q = (\pm1)^{n-1}\frac{\sin(\pi nK/N)}{\sin(\pi K/N)},\label{eq:sine}
\end{equation}
a form that is useful for computations.

The normalised monomials also generate a finite non-commutative torus. 
Let $a,b,c,d$ be integers, and define 
\begin{equation}X=E^{(a,b)},\quad Y=E^{(c,d)}.\label{eq:newgeometry}
\end{equation}
These obey the relation
\begin{equation}XY=\Q^{ad-bc}\,YX
\end{equation}
Moreover, the normalised monomials in $X$ and $Y$ are just examples of normalised monomials $E$ of $U$ and $V$ in the obvious way, i.e.,
putting $\X=\Q^{ad-bc}$ with square root $\X^{1/2}=\Q^{(ad-bc)/2}$, then
\begin{equation}\X^{-kl/2}X^kY^l=E^{(ka+lc,kb+ld)}.
\end{equation}
\bigskip

\paragraph{Clock and shift algebra}
The definition can be applied to the clock and shift operators $C$ and $S$ with a choice of $\q^{1/2}$. If the order $N$ of $\q$ is even, either square root can be chosen and $\q^{N/2}=-1$ for either choice. For $N$ odd, it is convenient to choose the root such that $\q^{N/2}=+1$. 

The normalised monomials of clock and shift are denoted $\e^{(m,n)}$ to prevent any confusion with more general cases.
The specific feature of the clock and shift generators is that $C^N=S^N=1$. This is preserved by the normalised monomials in the following sense.

\begin{lemma} Let $k$ be an integer so that $kn=0$ (mod $N$) and $km=0$ (mod $N$). Then
\begin{equation}(\e^{(m,n)})^k=1
\end{equation} 
\end{lemma}
\begin{proof}Using the periodicity,  $(\e^{(m,n)})^k=\e^{(km,kn)}=\e^{(0,0)}=1$. 
Note that the middle equality relies on the choice $\q^{N/2}=1$ in the case of $N$ odd.
\end{proof}

The monomials are periodic in $n$ and $m$; for $N$ odd  
\begin{equation} \e^{(m+N,n)}= \e^{(m,n)}, \quad \e^{(m,n+N)}= \e^{(m,n)},\label{eq:oddperiod}
\end{equation} 
while for $N$ even
\begin{equation} \e^{(m+N,n)}=(-1)^n \e^{(m,n)}, \quad \e^{(m,n+N)}=(-1)^m \e^{(m,n)}.\label{eq:evenperiod}
\end{equation}

\subsection{Transformations}
Two separate types of transformations of a torus algebra are defined here. Their geometric interpretation is discussed in Section \ref{sec:corr}.

\paragraph{Translations}
 Let $U$, $V$ generate a torus algebra with $\Q$ of order $N$.
\begin{definition} \label{def:translation}For each pair of integers $(j,n)$, the unitary operator $P_{(j,n)}=V^{-j}U^n$ determines the inner automorphism
\begin{equation} a\mapsto P_{(j,n)}^{\vphantom{-1}}\,a\,P_{(j,n)}^{-1}\label{eq:translation}\end{equation}
for all $a\in \langle U,V\rangle$.\end{definition}

In particular, $P_{(j,n)}^{\vphantom{-1}}UP_{(j,n)}^{-1}=\Q^jU$ and $P_{(j,n)}^{\vphantom{-1}}VP_{(j,n)}^{-1}=\Q^nV$. These transformations, for all $j$ and $n$, define an action of the abelian group $\Z_N\times\Z_N$ on  $\langle U,V\rangle$.
(Here $\Z_N\equiv\Z/N\Z$.) It is worth noting that the $P_{(j,n)}$ considered as operators in $\mathfrak h$ do not commute. They determine a representation of a central extension of $\Z_N\times\Z_N$ that is a type of Heisenberg group.

\paragraph{Modular transformations}
These are defined for the specific case of the algebra $M_N(\C)$ generated by $C$ and $S$.

\begin{definition} Let $A=\begin{pmatrix}a&b\\ c&d\end{pmatrix}\in\SL(2,\Z)$. The matrix determines an automorphism of $M_N(\C)$ by 
\begin{equation}T_A\colon C\mapsto \e^{(a,b)}, \quad S\mapsto \e^{(c,d)}\label{eq:fuzzytransform}
\end{equation}
This is called a modular transformation.
\end{definition}

This definition is similar to that in \cite{SS}.  However, it is important to note that a different normalisation factor has been used and the choice of this factor relates to whether $ N $ is even or odd. This is essential as it ensures that the spectrum of the scalar Laplace and Dirac operators on the fuzzy torus converge.  

A calculation shows that composing these automorphisms gives
\begin{equation}T_{A'}\circ T_A= T_{AA'}.
\end{equation} 
These transformations can be implemented by a unitary operator acting in $\mathfrak h$ according to \eqref{eq:trans}. This is because all representations of the clock and shift algebra are equivalent, according to Lemma \ref{lem:basis}. The next two examples give explicit formulas for two matrices that generate $\SL(2,\Z)$.

\begin{example} Let $A=\begin{pmatrix}1&1\\ 0&1\end{pmatrix}$. This is the transformation
\begin{equation}P_A^{\vphantom{-1}}CP_A^{-1}=\q^{-1/2}CS, \quad P^{\vphantom{-1}}_ASP_A^{-1}=S.
\end{equation}
An explicit formula for $P_A$ is
\begin{equation} P_A=\frac1{\sqrt N}\,  \sum_{n=0}^{N-1} \q^{n^2/2}\,S^n,
\end{equation}
and the matrix elements of $P_A=[(P_A)_{jk}]$, with indices $j,k$ from $0$ to $N-1$, are 
\begin{equation}(P_A)_{jk}=\frac1{\sqrt N} \q^{(j-k)^2/2}.
\end{equation}
\end{example}

\begin{example} Let $B=\begin{pmatrix}0&-1\\ 1&0\end{pmatrix}$. This is the transformation
\begin{equation}P_B^{\vphantom{-1}}CP_B^{-1}=S^{-1}, \quad P_B^{\vphantom{-1}}SP_B^{-1}=C.
\end{equation}
In this case, the matrix elements of $P_B$ are
\begin{equation} (P_B)_{jk}=\frac1{\sqrt N} \q^{jk},
\end{equation}
which is a Vandermonde matrix.
\end{example}


\subsection{Correspondence}\label{sec:corr}
This section gives some comments on the correspondence between the commutative torus and the finite non-commutative torus. The purpose is to guide geometrical thinking about the non-commutative torus; the ideas about a commutative limit are mostly heuristic at this stage.

The first point is that the monomials $E^{(m,n)}$ are non-commutative analogues of the functions $e^{i( m\theta+ n\phi)}$. In the commutative case $(m,n)$ are arbitrary integers and so the momentum space formed by the exponentials is the group $\Z\times\Z$. However, for the simplest case of the non-commutative torus in which $U$ and $V$ are clock and shift matrices, 
the momentum space has period $N$ according to \eqref{eq:oddperiod} and \eqref{eq:evenperiod}, and so is a discrete torus. Considering the $\e^{(m,n)}$ projectively, i.e., ignoring the phase factor, the multiplication law \eqref{eq:Emult} forms a group isomorphic to $\Z_N\times \Z_N$. 

Position space is dual to momentum space; the group dual to $\Z_N\times\Z_N$ is again $\Z_N\times\Z_N$. This suggests the finite non-commutative torus is analogous to a toroidal lattice. The translations on this lattice are the operators \eqref{eq:translation}, and the eigenvectors for these operators (commutatively the plane waves) are the monomials $\e^{(m,n)}$.

The action of the modular group is very similar in the commutative and non-commutative cases, the equation \eqref{eq:fuzzytransform} being the analogue of \eqref{eq:transform}. The difference is that for the clock and shift algebra the modular transformations act on a toroidal lattice.

The quantum integers have the following limit as $\Q\to1$. Suppose $\Q^{1/2}\to\sigma=\pm1$. Then for a fixed $n$ 
\begin{equation}[n]_\Q\to \sigma^{n-1} n.\label{eq:qilimit}
\end{equation}
Thus the $\Q^{1/2}\to 1$ limit of \eqref{eq:fuzzycommutators} is
\begin{equation}(ml-nk)\;E^{(m+k,n+l)}
\end{equation}
which is the analogue of the Poisson bracket $\Omega_0$, according to \eqref{eq:poissonexp}.

This suggests that one should think of an operator
\begin{equation}v_A\colon B\mapsto\frac1{\X^{1/2}-\X^{-1/2}} [A,B]
\end{equation}
as an analogue of the vector field $\Omega(A,\cdot)$ associated to the Hamiltonian $A$. It satisfies a Leibnitz identity,
\begin{equation}v_A(BC)=Bv_A C+ (v_A B)C.
\end{equation}
In the case considered in \eqref{eq:fuzzycommutators}, $\X=\Q$ and so the relevant Poisson structure is $\Omega_0$. In general, $\X=\Q^r$ and \eqref{eq:fuzzycommutators} generalises to
\begin{equation}\begin{aligned}\frac1{\X^{1/2}-\X^{-1/2}} [E^{(m,n)},E^{(k,l)}]
&=\frac{[(ml-nk)]_\Q}{[r]_\Q}\;E^{(m+k,n+l)}\\
\end{aligned}\label{eq:Qfuzzycommutators}\end{equation}

It is useful in places to extend Definition \eqref{eq:quantuminteger} of a quantum integer to allow a rational number $a/b$ as argument in place of $n$. This can be done as long as the definition of $\Q^{1/2b}$ is given. If  $\Q^{1/2b}=\pm e^{\pi iK/bN}$ then the sine formula \eqref{eq:sine}
again holds.

In \eqref{eq:Qfuzzycommutators}  one can write
\begin{equation}\frac{[(ml-nk)]_\Q}{[r]_\Q}=\left[\frac{ml-nk}r\right]_{\X},\end{equation}
 where the required root of $\X$ is defined by $\X^{1/2r}=\Q^{1/2}$.

\section{The fuzzy torus}\label{sec:fuzzy}

The fuzzy torus is a finite non-commutative torus with a Hilbert space in which $U$ and $V$ act on both the left and the right. This is achieved by adding a real structure to the non-commutative torus.

\begin{definition} \label{def:ncreal} A real structure for a finite non-commutative torus $(U,V,\mathfrak h)$ is an antiunitary map  $J\colon\mathfrak h\to\mathfrak h$ that is an involution, $J^2=1$, and obeys
\begin{equation}[JUJ,U]=0,\;[JUJ,V]=0, \;[JVJ,U]=0, \;[JVJ,V]=0.
\end{equation}
\end{definition}

The real structure in Definition \ref{def:ncreal} determines a right action of $U$ and $V$, written for $\psi\in\mathfrak h$ as
\begin{equation} \psi U=JU^*J\psi, \quad \psi V=JV^*J\psi.
\end{equation}
This commutes with the left action, making $\mathfrak h$ a bimodule over the algebra $\langle U,V\rangle$.
It also implies that commutators and anticommutators with elements of $\mathfrak h$ make sense, e.g.,
\begin{equation} [U,\psi]=U\psi-\psi U,\quad \{U,\psi\}=U\psi+\psi U.
\end{equation}
The right action is also unitary, i.e.,
\begin{equation}\psi U^*=\psi U^{-1},\quad \psi V^*=\psi V^{-1}.
\end{equation}
\begin{definition} A fuzzy torus is a finite non-commutative torus with a real structure.
\end{definition}

The Hilbert space of the fuzzy torus is the non-commutative analogue of the space of sections of a bundle over the torus. 
In the simplest case of a trivial line bundle the sections are just the complex functions and the real structure is complex conjugation. The fuzzy analogue of this is the following set of examples.
\begin{example} Let $(U,V,\mathfrak h')$ be a finite non-commutative torus. The Hilbert space of the fuzzy torus is
 $\mathfrak h=\End(\mathfrak h')$, with inner product $(\psi,\phi)=\tr(\psi^*\phi)$. Then $U, V\in \End(\mathfrak h')$ act by left and right multiplication, and the real structure is the Hermitian conjugate, $J=*$. These examples include the cases where $\mathfrak h'=\C^N$, $\mathfrak h=M_N(\C)$, and the multiplications are all matrix multiplication.
\end{example}

\subsection{The scalar Laplace operator}\label{sec:scalarlaplace}
The Laplace operator illustrates many of the features of the fuzzy torus in a simpler setting than the Dirac operator. One can define an analogue of the Laplace operator for the fuzzy torus and check that its eigenvalues converge to the eigenvalues of the corresponding differential operator on the flat torus. The modular transformations are also defined for the fuzzy torus and respect the limit.

The Laplace operator described here is the analogue of the differential Laplace operator defined on functions, or slightly more generally, sections of a vector bundle associated to a principal bundle with a discrete group. This is called the scalar Laplace operator. The metric on the torus is determined by two algebra elements $ X,Y \in C(T^2) $ according to \eqref{eq:cood-defs-gen}. Analogously, the scalar Laplace operator on a fuzzy torus is determined by $X, Y\in \langle U,V\rangle $ satisfying $ XY = \X YX $ with $ \X \in \C $.    
Let $\mathfrak h$ be the bimodule of a fuzzy torus determined by $U$ and $V$.
\begin{definition}
The  scalar Laplace operator on a fuzzy torus is the operator $\mathfrak h\to\mathfrak h$ given by
\begin{equation} \Delta_{X,Y}=\frac{-1}{(\X^{1/2}-\X^{-1/2})^2}\,\left(\,[X,[X^*,{\cdot}\,]]+[Y,[Y^*,\cdot\,]]\,\right),
\end{equation}
with $ X, Y \in \langle U, V \rangle $.
\end{definition}
This operator is self-adjoint and so has real eigenvalues.
The construction is respected by transformations that act on $\mathfrak h$. For example, let $P\in\langle U,V\rangle$ such that $PXP^{-1}=X'$, $PYP^{-1}=Y'$. Then $P$ acts in $\mathfrak h$ by conjugation, and $(P\cdot P^{-1})\Delta_{X,Y}(P^{-1}\cdot P)=\Delta_{X',Y'}$.

\begin{example}\label{ex:CS}
In the simple case of the algebra of Example \ref{ex:cs}, with $\mathfrak h=M_N(\C)$, $J=*$, $X=C$, $Y=S$, the eigenvalues can be computed using \eqref{eq:fuzzycommutators},
\begin{equation} \Delta_{C,S}\,\e^{(k,l)}=([k]_\q^2+[l]_\q^2)\,\e^{(k,l)}.\label{eq:ev-fuzzysquare}
\end{equation}
In the limit  $\q\to1$ with $k$ and $l$ fixed,
\begin{equation}[k]_\q^2+[l]_\q^2\to k^2+l^2.\label{eq:naivelimit}
\end{equation}
These are the eigenvalues for the  Laplacian
\begin{equation}\Delta_\delta=-\left(\pd{\mathstrut^2}{\phi^2}+\pd{\mathstrut^2}{\theta^2}\right)
\end{equation}
acting on functions on the square flat torus with metric $\delta$, the eigenvalue equation being
\begin{equation} \Delta_\delta\; e^{i(k\theta+l\phi)}=(k^2+l^2)\, e^{i(k\theta+l\phi)}.\label{eq:ev-square}
\end{equation}

Looking at the spectrum as a whole, the picture is a little more complicated. Let $\q=e^{2\pi i K/N}$, with $0<K<N/2$ and coprime to $N$, so that $\q$ has order $N$. The sine curves \eqref{eq:sine} would each have $K$ zeroes if the argument were continuous; for the integer values therefore the eigenvalues are close to the value $0$ for  $K^2$ times as $k$ and $l$ cycle through their $N$ values. However, the local minima are not all the same. 

For example, if $N$ is odd and $K=2$ then   $\q^{1/2}=e^{2i\pi/N}$ and
\begin{equation} [0]_\q^2=0\quad\text{ but }\quad \left[\frac {N\pm 1}2\right]_\q^2= \left[\frac  12\right]_\q^2.
\end{equation}
Shifting the indices $k$ and $l$ by defining integers $k'=k-N/2-1/2$ and  $l'=l-N/2-1/2$ gives
\begin{equation}
[k]^2_\q=\left[k'+\frac12\right]^2_\q, \quad [l]^2_\q=\left[l'+\frac12\right]^2_\q.
\end{equation}
Now taking a limit $\q^{1/2}\to1$ with constant $(k',l)$, $(k,l')$ or $(k',l')$ leads to the three other limits
\begin{equation} \left(k'+\frac12\right)^2+l^2,\quad k^2+\left(l'+\frac12\right)^2,\quad\text{or}\quad  \left(k'+\frac12\right)^2+\left(l'+\frac12\right)^2.
\end{equation} 
These spectra are characteristic of the Laplacian on line bundles over the torus, as will be explained in Section \ref{sec:linebundles}. Note that this analysis does not work if $K$ is odd and fixed in the limit; in those cases $\q^{1/2}\to-1$ and $[1/2]_\q$ does not converge at all.
\end{example}
\bigskip

The limits in the previous example can be understood using the heuristics for the commutative limit from Section \ref{sec:corr}. For the eigenvectors $C\rightsquigarrow e^{i\theta}$, $S\rightsquigarrow e^{i\phi}$, and for the scaled
commutator $v_X\rightsquigarrow\Omega_0(X,\cdot\,)$. Then, assuming  $\q^{1/2}\to1$,
\begin{equation}\begin{gathered} \frac1{\q^{1/2}-\q^{-1/2}}[C,\cdot\,]\rightsquigarrow-ie^{i\theta}\pd{}\phi,\quad   \frac1{\q^{1/2}-\q^{-1/2}}[C^*,\cdot]\rightsquigarrow ie^{-i\theta}\pd{}\phi\\
 \frac1{\q^{1/2}-\q^{-1/2}}[S,\cdot\,]\rightsquigarrow ie^{i\phi}\pd{}\theta,\quad   \frac1{\q^{1/2}-\q^{-1/2}}[S^*,\cdot]\rightsquigarrow -ie^{-i\phi}\pd{}\theta.\\
\end{gathered}\label{eq:commutatorlimit}
\end{equation}
With these replacements, equation \eqref{eq:ev-fuzzysquare} becomes \eqref{eq:ev-square}.
\bigskip

\begin{example}\label{ex:genXY}
More generally, with the same Hilbert space and algebra as Example \ref{ex:CS}, but defining the geometry with $X=\e^{(a,b)}$, $Y=\e^{(c,d)}$, so that $ \X = \q^{ad-bc}$,
one can check that this geometry is a fuzzy version of \eqref{eq:transform}. 

The eigenvectors and eigenvalues of $ \Delta_{X,Y} $ are
\begin{equation} \Delta_{X,Y}\,\e^{(k,l)}=\frac{[al-bk]_\q^2+[dk-cl]_\q^2}{[ad-bc]^2_\q}\,\e^{(k,l)} \label{eq:laplaceev-integral}
\end{equation}
for integers $(k,l)\in\Z\times\Z$, the eigenspaces being periodic in $k$ and $l$ with period $N$.

Using the same technique as for the square torus example, the limiting eigenvalue equation is
\begin{equation} \Delta_g\; e^{i(k\theta+l\phi)}=\frac{(al-bk)^2+(dk-cl)^2}{(ad-bc)^2}
\, e^{i(k\theta+l\phi)}.\label{eq:ev-integral}
\end{equation}
with the Laplacian
\begin{equation} \Delta_g=\frac{-1}{(ad-bc)^2}\left((b^2+d^2)\pd{\mathstrut^2}{\theta^2}-2(ab+cd)\pd{\mathstrut^2}{\theta\partial\phi}+(a^2+c^2)\pd{\mathstrut^2}{\phi^2}\right).\label{eq:laplacian}
\end{equation}
Writing $\Delta_g=-g^{ab}\partial_a\partial_b$, the coefficients $g^{ab}$ are exactly the inverse of the integral metric \eqref{eq:cood-defs-gen}.
\end{example}

In some further examples the algebra $ \langle U, V \rangle $ may be a proper subalgebra of $M_{N}(\C)$. In these cases it is possible to decompose $ \mathfrak{h} $ into subspaces that are still bimodules over the algebra.

\begin{example} \label{ex:C2S}  This uses the algebra of Example \ref{ex:Csquard}, with $\mathfrak h = M_N(\mathbb{C})$, $J=*$, $X=C^2$, and $Y=S$. The eigenvalue of eigenvector $C^kS^l$ is 
\begin{equation}\lambda_{k,l}=\frac{[k]_\q^2+[2l]_\q^2}{[2]^2_\q}=\left[\frac k2\right]_\Q^2+[l]_\Q^2.
\end{equation}

\begin{itemize}
\item For $N$ odd, $\{[2l]^2_\q \mid l=1,\ldots, N\}=\{[l]^2_\q\mid l=1,\ldots, N\}$. Thus the spectrum is, counterintuitively, the same as the unit square torus, except for the overall normalisation factor of $[2]^2_\q$. 

\item For $N$ even, $\Q$ has order $N/2$ and so $[l+N/2]^2_\Q=[l]^2_\Q$, which means that the $l$ term in the eigenvalue formula ranges over the set  $\{[l]^2_\Q \mid l=1,\ldots, N/2\}$, with multiplicity two. 

The torus algebra is $\langle U,V\rangle\cong M_{N/2}(\C)\oplus  M_{N/2}(\C)$, a proper subalgebra of $M_N(\C)$.  The Hilbert space splits as a bimodule as $\mathfrak h=\mathfrak h_0\oplus\mathfrak h_1$ with $\mathfrak h_0$ spanned by the monomials with even powers of $C$ and $\mathfrak h_1$ the odd powers. Then the eigenvalues on $\mathfrak h_0$ are
\begin{equation}\lambda_{2m,l}=[m]^2_\Q+[l]_\Q^2\end{equation}
for $m=1,\ldots,N/2$. The eigenvalues on $\mathfrak h_1$ are
\begin{equation}\lambda_{2m-1,l}=\left[m-\frac12\right]^2_\Q+[l]_\Q^2.\end{equation}
\end{itemize}
\end{example}

\begin{example} \label{ex:C2S2} Consider the fuzzy torus defined by $\mathfrak h=M_N(\C)$ with $N$ even and $J=*$, $X=U=C^2$, $Y=V=S^2$. 
Then $ \X = \Q = \q^4$, which has order $N/2$ (if $N/2$ is odd) or $N/4$ (if $N/2$ is even). The eigenvalue of eigenvector $C^kS^l$ is 
\begin{equation}\lambda_{k,l}=\left[\frac k2\right]_\Q^2+\left[\frac l2\right]_\Q^2. \label{eq:Laplacesq}\end{equation}

The Hilbert space splits into four subspaces with even or odd powers of $C$ and $S$,
\begin{equation}\mathfrak h=\mathfrak h_{00}\oplus\mathfrak h_{01}\oplus\mathfrak h_{10}\oplus\mathfrak h_{11},\end{equation}
with $\mathfrak h_{hj}$ spanned by monomials $C^kS^l$ satisfying $(k,l)=(h,j) \mod 2$. Note that these subspaces can also be characterized as the eigenspaces of the adjoint action of $C^{N/2}$ and $S^{N/2}$. The eigenvalues on these subspaces are
\begin{equation}\begin{aligned}\mathfrak h_{00}\colon\quad&\lambda_{2m,2n}=\left[m\right]_\Q^2+\left[n\right]_\Q^2\\
\mathfrak h_{01}\colon\quad&\lambda_{2m-1,2n}=\left[m-\frac12\right]_\Q^2+\left[n\right]_\Q^2\\
\mathfrak h_{10}\colon\quad&\lambda_{2m,2n-1}=\left[m\right]_\Q^2+\left[n-\frac12\right]_\Q^2\\
\mathfrak h_{11}\colon\quad&\lambda_{2m-1,2n-1}=\left[m-\frac12\right]_\Q^2+\left[n-\frac12\right]_\Q^2.
\end{aligned}\label{eq:foursectors}\end{equation}
for integer $m,n$, the eigenspaces recurring with period $N/2$. The formulas for the eigenvalues each have period $N/4$, however.

The torus algebra $\langle C^2,S^2\rangle$ has dimension $N^2/4$. It can be characterised as the subalgebra of $M_N(\C)$ that commutes with $C^{N/2}$ and $S^{N/2}$. 

\begin{itemize}\item In the case where $N/2$ is odd, $C^{N/2}$ and $S^{N/2}$ anti-commute and form a Clifford algebra. Since $\Q$ has order $N/2$, $\langle C^2,S^2\rangle\cong M_{N/2}(\C)$. The left (or right) actions of $C^{N/2}$ and $S^{N/2}$ on $\mathfrak h$ permute the four subspaces in \eqref{eq:foursectors} and also commute with the Laplacian.  This shows that the four spectra are in fact the same, the apparently different formulas  being an instance of the phenomenon noted already in Example \ref{ex:CS}. Each bimodule $\mathfrak h_{hj}$ is isomorphic to an instance of Example \ref{ex:CS}, with data $N'=N/2$, $C'=C^2$, $S'=S^2$ and $\q'=\Q$. 

\item In the case where $N/2$ is even, $C^{N/2}$ and $S^{N/2}$ are central elements of $\mathcal A=\langle C^2,S^2\rangle$. The eigenvalues $\pm1$  of each provide the central idempotents $\pi_{hj}=\frac14(1+(-1)^hC^{N/2})(1+(-1)^jS^{N/2})$. Define the subalgebras $\mathcal A_{hj}=\pi_{hj}\mathcal A$. Then there is a splitting 
\begin{equation} \mathcal A=\mathcal A_{00}\oplus\mathcal A_{01}\oplus\mathcal A_{10}\oplus\mathcal A_{11}.
\end{equation}
There is a canonical isomorphism of $\mathcal A_{00}$ to the clock and shift algebra $M_{N/4}(\C)=\langle C',S'\rangle$ given by
$\pi_{00}C^{2k}S^{2l} \mapsto C'^kS'^l$. The other summands are also isomorphic to $M_{N/4}(\C)$ but not in a canonical way. For example, the map $a\mapsto S^naS^{-n}$ is an isomorphism $\mathcal A_{10}\to \mathcal A_{00}$ for any odd integer $n$.

The bimodule $\mathfrak h_{00}$ is isomorphic to the algebra acting on itself, but the bimodules $\mathfrak h_{hj}$ are not isomorphic to each other.
\end{itemize}
\end{example}

\subsection{Line bundles}\label{sec:linebundles}
Examples \ref {ex:C2S} and \ref{ex:C2S2} consider splitting the Hilbert space spanned by clock and shift operators into subspaces with odd and even powers of these operators. Here, the geometry of the corresponding algebras and bimodules in the commutative case is explored and generalised. 

The functions $U=e^{2i\theta}$ and $V=e^{2i\phi}$ can be understood as pull-backs of $C=e^{i\theta}$ and $S=e^{i\phi}$ by the covering map $(\theta,\phi)\mapsto(2\theta,2\phi)$. The Hilbert space  of complex-valued functions on $T^2$ with metric $g$ is denoted $\mathfrak h=L^2(T^2,g)$. 
This Hilbert space also splits into subspaces $\mathfrak h_{hj}$, with $h,j\in\{0,1\}$, spanned by odd or even powers of $e^{i\theta}$ and $e^{i\phi}$. A function $\psi\in \mathfrak h_{hj}$ is the pull-back of a section $\Psi$ of a complex line bundle over $T^2$, i.e., $\psi=\Psi(2\theta,2\phi)$ with
\begin{equation} \quad \Psi(\theta+2\pi,\phi)=(-1)^h\Psi(\theta,\phi),\quad \Psi(\theta,\phi+2\pi)=(-1)^j\Psi(\theta,\phi).
\end{equation}
The line bundle is specified by the periodic or anti-periodic boundary conditions determined by $h$ and $j$.

\paragraph{General case}
Considering a more general case puts the about considerations in context and provides some results that are useful later.

A map of manifolds that is  a local isomorphism (i.e., an immersion of manifolds of the same dimension) and onto is called a covering map. Let $\mathcal M$ be a connected manifold and 
\begin{equation} \tau\colon \widehat {\mathcal M}\to \mathcal M
\end{equation}
a regular covering map, with $\widehat {\mathcal M}$ also connected. Here, `regular' means that $\tau_\bullet\bigl(\pi_1(\widehat {\mathcal M})\bigr)$ is a normal subgroup of $\pi_1(\mathcal M)$. Then, putting 
\begin{equation} \mathcal G=\pi_1(\mathcal M)/\tau_\bullet\bigl(\pi_1(\widehat {\mathcal M})\bigr),
\end{equation}
the covering is a principal $\mathcal G$-bundle, with a right action of $\mathcal G$ on $\widehat {\mathcal M}$.
If $\chi$ is a unitary character for $\mathcal G$ (a group homomorphism $\mathcal G\to \U(1)$) then the associated bundle construction determines a Hermitian line bundle $L_\chi$ over $\mathcal M$.  In this construction, a point of the bundle is an equivalence class $[(z,\zeta)]\subset \widehat {\mathcal M}\times\C$ under the relations $(z,\zeta)\sim (zg,\chi(g^{-1})\zeta)$ for all $g\in \mathcal G$.
The sections of $L_\chi$ correspond to functions $\psi$ on $\widehat {\mathcal M}$ such that
\begin{equation} \psi(zg)=\chi(g^{-1})\psi(z) \quad\text{for all } g\in \mathcal G.\label{eq:equivariant}
\end{equation}
Since $\mathcal G$ is a discrete group, the principal bundle $\widehat {\mathcal M}$ has a unique connection on it and therefore the line bundle $L_\chi$ has a uniquely-defined covariant derivative. 

Now suppose that $\mathcal G$ is a finite group and that $\mathcal M$ has a volume form $\upsilon$. Then $\widehat {\mathcal M}$ has the volume form $\widehat\upsilon=\frac 1 {|\mathcal G|}\tau^\bullet\upsilon$.  
Define $\mathfrak h=L^2(\widehat {\mathcal M},\widehat\upsilon)$, the Hilbert space of complex functions on the covering space. The group $\mathcal G$ has a unitary left action in this Hilbert space by $\psi\mapsto g\psi$ with
\begin{equation} (g\psi)(z)=\psi(zg).\end{equation}
Let  $\mathfrak h_\chi$ be the subspace in which $\mathcal G$ acts by the character $\chi$, i.e.,
\begin{equation} (g\psi)(z)=\chi(g)\psi(z)\label{eq:linesplit}.\end{equation}
Comparing with \eqref{eq:equivariant} shows that  $\mathfrak h_\chi$ is isomorphic to $L^2(\mathcal M,L_{\overline\chi},\upsilon)$, the sections of the line bundle determined by the conjugate character $\overline\chi(g)=\chi(g^{-1})$.

If $\mathcal G$ is an abelian group then all the irreducible representations of $\mathcal G$ are determined by unitary characters and so
\begin{equation} \mathfrak h=\bigoplus_{\chi\in \mathcal G^*}\mathfrak h_\chi, \label{eq:hdecomp}
\end{equation} 
where $\mathcal G^*$ is the dual of $\mathcal G$.

\section{The Dirac operator on a manifold}\label{sec:dirac}
This section presents the general formalism for Dirac operators on manifolds in a way that is amenable for the construction of non-commutative analogues. 

\subsection{Spin structures}\label{sec:spinstructures}
The definition of a spin structure is independent of the choice of metric, so this is a good place to start the discussion of spin geometry. Let $\mathcal M$ be an oriented manifold of dimension $m$. A frame at a point $x$ is a linear isomorphism $e_x\colon \R^m\to T\mathcal M_x$ preserving the orientation; the collection of all frames for all $x$ is a principal bundle $\pi\colon \mathcal E\to \mathcal M$ called the (oriented) frame bundle. It has structure group $\GL^+(m)$; an element $l$ of this group acting on frames by
\begin{equation}e_x\mapsto e_x\, l.
\end{equation}

The idea behind the definition of spin structures is to describe the possible liftings of $\mathcal E$ to principal bundles with a structure group that covers $\GL^+(m)$.  This will be described more fully in Section \ref{sec:diracmfd}. 

The inclusion $\SO(m)\subset \GL^+(m)$ is a homotopy equivalence and so \\ $\pi_1(\GL^+(m))=\pi_1(SO(m))$.  For  $m\ge3$ this group is $\Z_2$ and so the universal covering group of  $\GL^+(m)$ is a 2-1 covering. For $m<3$ this is not true, which complicates the definition of spin structure in these cases.

So first assume that $m\ge3$. A spin structure \cite{Kirby} is a cohomology class $s\in H^1(\mathcal E;\Z_2)$ having the property that the restriction to a fibre is non-trivial, i.e., if $c$ is the non-trivial 1-cycle, then $s(c)=1$. If spin structures exist, then $\mathcal M$ is called a spin manifold. The group $H^1(\mathcal M;\Z_2)$ acts freely on the set of spin structures by
\begin{equation}s\mapsto s+\pi^\bullet t\label{eq:spindiff}\end{equation}
for $t\in H^1(\mathcal M;\Z_2)$, and so the number of spin structures is the number of elements in this group. 

A standard result \cite{Hatcher} is that for any connected space $X$,
\begin{equation} H^1(X;\Z_2)\cong\Hom(H_1(X;\Z),\Z_2)\cong\Hom(\pi_1(X),\Z_2).\label{eq:stdresults}
\end{equation}
Therefore, a spin structure on $\mathcal E$ is equivalent to a homomorphism   $\pi_1(\mathcal E)\to\Z_2$.

This formalism can be extended by adding extra trivial dimensions to the tangent space. The tangent space is replaced by $T\mathcal M_x\oplus \R^k$, for some fixed $k$, and a frame is now an isomorphism 
\begin{equation}e_x\colon \R^{m+k}\to T\mathcal M_x\oplus \R^k\label{eq:frame}\end{equation} 
preserving an orientation on this augmented bundle. The principal bundle $\mathcal E$ is now the set of all frames of this sort. Examples are where $\mathcal M$ is immersed in an oriented manifold $\mathcal N$ of dimension $m+k$ with a given trivialisation of the normal bundle. Since the inclusion $\GL^+(m)\subset\GL^+(m+k)$ is a homotopy equivalence for $m\ge3$, this does not change the set of spin structures. 

The simplest example is to take $\mathcal N=\mathcal M\times\R^k$. This example allows the uniform treatment of manifolds of dimensions $0,1$ and $2$, by stabilising the tangent bundle so it has dimension $3$ or more, then applying the formalism already outlined. This defines spin structures for these manifolds.

Let $\tau\colon \mathcal M\to \mathcal M'$ be a covering map. It extends to a map of the frame bundles and so if $s'$ is a spin structure on $\mathcal M'$, its pull-back $s=\tau^\bullet s'$ is a spin structure on $\mathcal M$. 

\begin{example}\label{ex:torusss} The torus $T^2$ has the frame field $e_x(u,v)=u\pd{}\theta+v\pd{}\phi$. Since this is constant in $x=(\theta,\phi)$ it is called a parallel frame field. The frame field is a section of the frame bundle and $t=e^\bullet(s)\in H^1(T^2;\Z_2)$ parameterises the spin structures. If $c_1$ and $c_2$ are the 1-cycles along the two axes, then $\sigma=((t(c_1),t(c_2))\in \Z_2\times\Z_2$ is an explicit parameterisation of the four spin structures on the torus. The spin structure $s_0$ labelled by $\sigma=(0,0)$ is called the Lie spin structure and one can write
\begin{equation} s-s_0=\pi^\bullet t.
\end{equation}

 Let $A\in \GL^+(2,\R)\cap M_2(\Z)$. This determines a covering map $\widetilde A\colon T^2\to T^2$ and hence a pull-back spin structure $\widetilde A^\bullet s$. The key fact about the torus is that the Lie spin structure is invariant under this pull-back. This follows from the fact that parallel frame fields with the same orientation are related by a $\GL^+(2,\R)$ transformation and so can be deformed to each other. Thus
\begin{equation} \widetilde A^\bullet s-s_0=\widetilde A^\bullet \pi^\bullet t=\pi^\bullet \widetilde A^\bullet t
\end{equation}
and so spin structure $s'=\widetilde A^\bullet s$ is parameterised by $t'=\widetilde A^\bullet t$. 

Denote by $\lBracket n\rBracket_2$ the  integer $n$ modulo 2. Explicitly, 
\begin{multline} \sigma'=((t'(c_1),t'(c_2))=((t(A_\bullet c_1),t(A_\bullet c_2))=((t(a c_1+c c_2),t(b c_1+d c_2))\\
=(\lBracket a\rBracket_2\,t(c_1)+\lBracket c\rBracket_2\,t(c_2),\lBracket b\rBracket_2\,t(c_1)+\lBracket d\rBracket_2\,t(c_2))=\lBracket {A}^T\rBracket_2\,\sigma
\label{eq:sspb}\end{multline}
with $\lBracket A^T\rBracket _2$ the transpose of the matrix $A$ with entries modulo 2.

The three spin structures $(0,1)$, $(1,0)$ and $(1,1)$ are permuted by isomorphisms and all of them bound a spin three-manifold.
\end{example}

\subsection{The Dirac operator}\label{sec:diracmfd}
On a Riemannian manifold the frames of the last section can be specialised to orthonormal frames, and the groups reduced to orthogonal groups and their spin covering groups. This has the advantage that it allows the construction of vector bundles of spinors but at the cost of introducing a dependence on the choice of Riemannian metric.

The group $\SO(m)$ has a 2-1 covering by the homomorphism $\Phi\colon\Spin(m)\to\SO(m)$. For $m\ge 3$, $\Spin(m)$ is the universal connected cover. The group $\Spin(2)$ is defined to be $U(1)$ ($\cong\SO(2)$) with $\Phi$ the 2-1 cover, while $\Spin(1)$ and $\Spin(0)$ are both $\Z_2$. These definitions are compatible with the stabilisation described in the last section. While $\SO(2)$ also has $r$-fold covers for  $r>2$, the corresponding $r$-spin structures are not considered here.

Now suppose that $\mathcal M$ is a Riemannian manifold with an orientation.
There is a sub-bundle $\mathcal O$ of $\mathcal E$ consisting of the oriented orthonormal frames. This bundle has structure group $\SO(m)$. A spin structure $s$ on $\mathcal E$ restricts to an element of $H^1(\mathcal O;\Z_2)$ and determines a homomorphism $\pi_1(\mathcal O)\to\Z_2$ using \eqref{eq:stdresults}. This can be used to construct a two-fold covering $\eta\colon \mathcal F\to \mathcal O$ and hence a principal bundle $\pi\circ\eta\colon \mathcal F\to \mathcal M$ with structure group $\Spin(m)$, which is a lifting of the principal $\SO(m)$ bundle $\mathcal O$ \cite{Barrett:1991aj}. Note that the covering $\eta$ contains the information about the spin structure \cite{Kirby}.

A bundle of spinors is a Hermitian vector bundle (a complex vector bundle with a Hermitian inner product on each fibre) $\mathcal W$ associated to $\mathcal F$. This is constructed from the vector space $\C^p$ with the standard Hermitian inner product, on which the $m$-dimensional Clifford algebra acts irreducibly. This is specified by
$\{\gamma^i, i=1,\ldots,m\}$, a fixed set of anti-Hermitian $p\times p$ matrices, called gamma matrices, which satisfy
\begin{equation}\gamma^i\gamma^j+\gamma^j\gamma^i=-2\delta^{ij}.
\end{equation}
For even $m$, $p=2^{m/2}$, while for odd $m$, $p=2^{(m-1)/2}$.

The gamma matrices determine a unitary action of $\Spin(m)$ on $\C^p$ which is used to construct the spin bundle according to the associated bundle construction.
The action of a spin transformation $Z$ defines the corresponding rotation $R=\Phi(Z)\in\SO(m)$ according to
\begin{equation}\sum_j(R^{-1})^{ij}\gamma^j=Z\gamma^iZ^{-1}\label{eq:intertwine}\end{equation}
with $(R^{-1})^{ij}$ the matrix components of $R^{-1}$. 

The Clifford module $\C^p$ has some further properties. 
It has a real structure, which is an antiunitary map $j\colon \C^p\to\C^p$ commuting with the action of  $\Spin(m)$, with $j^2=\epsilon=\pm1$.

For even $m$, it also has a chirality operator $\gamma\colon\C^p\to\C^p$ which is linear, unitary, satisfies $\gamma^2=1$, and commutes with the action of   $\Spin(m)$. It relates to the real structure by $j\gamma j^{-1}\gamma =\epsilon''=\pm1$.  The signs are determined by the dimension of the spin group $m$, mod 8. For even dimensions they are given by the following table; for odd dimensions see \cite{Barrett:2015naa}.
\begin{equation}\label{eq:table}
\begin{array}{c|c|c|c|c}
m&0&2&4&6\\
\hline
\epsilon&1&-1&-1&1\\
\epsilon''&1&-1&1&-1
\end{array}
\end{equation} 
The Clifford multiplication on $\C^p$ is the map 
\begin{equation}\begin{aligned}\rho_0\colon (\R^{m})^*&\to\End(\C^p)\\
\rho_0(\nu)&=\sum_i\nu_i\,\gamma^i \end{aligned}\label{eq:cliffmodule}
\end{equation}
with $\nu_i\in\R$ the components of $\nu$.

These structures on $\C^p$ carry over to similar properties of the spin bundle: there is a real structure $j$ and chirality operator $\gamma$ on each fibre $\mathcal W_x$. A Clifford multiplication is the action of cotangent vectors on $\mathcal W_x$,
\begin{equation}\rho_x\colon TM_x^*\to\End(\mathcal W_x),\end{equation}
satisfying at each point $x$ the relations 
\begin{gather}\rho_x(\nu)\rho_x(\nu')+\rho_x(\nu')\rho_x(\nu)=-2 g(\nu,\nu')\\
j\rho_x(\nu)=\rho_x(\nu)j
\end{gather}
for all $\nu$ and $\nu'$.

A point  $f\in\mathcal F_x$ determines a linear isomorphism $f\colon\C^p\to \mathcal W_x$ called a spin frame. Explicitly, this is obtained by mapping $\zeta\in\C^p$ to the equivalence class $[(f,\zeta)]$ that defines a point of $\mathcal W_x$. This linear map determines the point $f$ uniquely and so $\mathcal F$ is often called the bundle of spin frames.

Denote the corresponding orthonormal frame $e=\eta(f)$. The Clifford multiplication on the spinor bundle at point $x$ is defined by 
\begin{equation}\rho_x(\nu)=f\,\rho_0(e^*\nu) f^{-1}.\label{eq:clifford}\end{equation}

\begin{lemma}\label{lem:cliff}
The Clifford multiplication is independent of the choice of frame.
\end{lemma} 
\begin{proof}
The intertwining property \eqref{eq:intertwine} can be written 
\begin{equation}\rho_0\left(\nu\right) =Z\rho_0\left(R^*\nu\right) Z^{-1}.\label{eq:intertwineZ}
\end{equation}
If two spin frames differ by $f'=fZ$, and hence $e'=eR$, then this implies
\begin{equation}f\,\rho_0(e^*\nu) f^{-1}=f'\,\rho_0(e'^*\nu) f'^{-1}.
\end{equation}
and so $\rho_x$ is independent of the choice of spin frame.
\end{proof}

Denote the smooth sections of a bundle by $\Sigma$ and define the multiplication map $\mu\colon\Sigma(T^*\mathcal M\otimes  \mathcal W)\to \Sigma(\mathcal W)$  by 
\begin{equation}\left(\mu(\nu\otimes w)\right)_x=\rho_x(\nu_x)w_x.\label{eq:mu}\end{equation}

A connection on the principal bundle $\mathcal O$ determines a unique connection on $\mathcal F$ and hence a covariant derivative on sections of $\mathcal W$. It will be assumed that $\mathcal O$ has a connection. For the torus this will be the Levi-Civita connection, which is uniquely determined by the metric. 

If $\nabla\colon  \Sigma(\mathcal W)\to \Sigma(T^*\mathcal M\otimes  \mathcal W)$ is a covariant derivative on $\mathcal W$, then the Dirac operator on  $\Sigma(\mathcal W)$ is defined by
\begin{equation}D\psi=\mu\nabla\psi.\label{eq:commDirac}
\end{equation}

This formalism can be extended to the case where the tangent bundle is augmented with extra trivial directions $\R^k$ and a spinor bundle $\mathcal W$ that accommodates the Clifford multiplication for the extra directions.
 The Clifford multiplication on this bundle is 
\begin{equation}\check\rho_x\colon (T\mathcal M_x\oplus \R^k)^*\to\End( \mathcal W_x),\label{eq:extendedClifford}\end{equation}
and it is assumed that the extra directions have the standard metric and orientation, and are orthogonal to the tangent space. Then if  $P_x\colon T\mathcal M_x\oplus \R^k\to T\mathcal M_x$ is the orthogonal projection, putting
\begin{equation}\rho_x=\check\rho_x\, P_x^*\label{eq:projcliff}
\end{equation}
defines a Dirac operator again, using \eqref{eq:commDirac}. The connection is again unique if one assumes that constant vectors in the extra directions are covariantly constant. Essentially, the extra directions act on the (possibly larger) spinors but play no role in the Dirac operator.

Formula \eqref{eq:clifford} generalises to this case, giving
\begin{equation}\check\rho_x(\omega)=f\,\rho_0(e^*\omega) f^{-1}\label{eq:gencliff}\end{equation}
for the Clifford multiplication of $\omega\in (T\mathcal M_x\oplus\R^k)^*$, using a frame $e$ for the $\SO(m+k)$ bundle at point $x$, a frame $f$ for the $\Spin(m+k)$ bundle, and $m+k$ gamma matrices.

\subsection{Frame fields}\label{sec:framefields}
The Dirac operator can be written using a trivialisation of the orthonormal frame bundle $\mathcal O$, i.e., a choice of an orthonormal frame $e_x$ for each point  $x$ of (a subset of) $\mathcal M$. This is called a frame field. 

Adding trivial directions to the tangent bundle gives a greater flexibility in the formulas that is important in this paper. This is because the frames do not have to respect the direct sum decomposition of $T\mathcal M_x\oplus \R^k$.

Let $\xi_i$ be the standard basis vectors of $\R^{m+k}$. Then the orthogonal frame $e$ determines $m+k$ vector fields on $\mathcal M$, denoted $v_i$, by projecting onto the part tangent to $\mathcal M$. The definition is
\begin{equation}(v_i)_x=P_xe_x(\xi_i).\label{eq:vi}
\end{equation}
The inverse metric tensor is given by the expression
\begin{equation}\sum_{i=1}^{m+k} v_i \otimes v_i\in T\mathcal M\otimes T\mathcal M.\label{eq:inversemetric}
\end{equation}
The Dirac operator can be expressed in terms of the frame fields.

\begin{lemma} Let $\mathcal U\subset \mathcal M$ be an open subset on which there is a frame field $e\colon \mathcal U\to \mathcal O$ covered by a spin frame field $f\colon \mathcal U\to \mathcal F$. Then the Dirac operator on sections of the spinor bundle restricted to  $\mathcal U$ is
\begin{equation}D\psi=\sum_{i=1}^{m+k}f \gamma^i f^{-1}\,\nabla_{v_i}\psi.\label{eq:diraca}
\end{equation}\label{lem:localDirac}
\end{lemma}

\begin{proof} The conclusion follows immediately from the formula
\begin{equation}\mu=\sum_iv_i\otimes f\gamma^if^{-1}.
\end{equation}
To prove this formula, suppose that $\nu\in \Sigma(T^*\mathcal M)$ and $w\in\Sigma(\mathcal W)$. Then
\begin{equation}\left(\sum_iv_i\otimes f\gamma^if^{-1}\right)(\nu\otimes w)=\sum_i\nu(v_i)f\gamma^if^{-1}w.
\end{equation}
At the point $x$ this spinor is
\begin{multline}\sum_i\nu_x(P_xe_x\xi_i)\,f_x\gamma^if_x^{-1}w_x=\sum_i(e_x^*P_x^*\nu_x)(\xi_i)\,f_x\gamma^if_x^{-1}w_x=\\
f_x\rho_0(e_x^*P_x^*\nu_x)f_x^{-1}w_x=\check\rho(P_x^*\nu_x)w_x=\rho(\nu_x)w_x=\mu(\nu\otimes w)_x,
\end{multline}
using  \eqref{eq:vi}, \eqref{eq:cliffmodule}, \eqref{eq:gencliff}, \eqref{eq:projcliff},  and  \eqref{eq:mu}.
\end{proof}

\paragraph {Global frame fields}
From here on it will be assumed that the frame field $e$ is global, i.e., exists smoothly over the whole of the manifold $\mathcal M$. This is appropriate for the torus, and looks to be a useful geometric starting point for the construction of fuzzy spaces more generally.

Suppose $\mathcal M$ is a Riemannian manifold, with its tangent space possibly augmented by trivial directions, as above.
An orthonormal frame field $e$ for $T\mathcal M\oplus \R^k$ determines a particular spin structure $s_c$ called the canonical spin structure. This is the spin structure such that $e^\bullet(s_c)=0$. 
The principal bundle of spin frames is the product bundle $\mathcal F=\mathcal M\times \Spin(m+k)$ and the two spin frames projecting to $e$ are defined to be $f=\pm 1$, with $1$ standing for the identity element of $ \Spin(m+k)$. The associated spinor bundle is as follows.

\begin{definition}\label{def:trivial} Let $e$ be an orthonormal frame field for $T\mathcal M\oplus \R^k$. Then the trivial spinor bundle for the frame field $e$ is defined to be the product vector bundle $\mathcal W^0=\mathcal M\times \C^p$. The Clifford multiplication on this bundle is defined according to \eqref{eq:gencliff} with $f=1$.
\end{definition}

Using the same frame field, one can describe all of the other spin structures on $\mathcal M$. 
Let $s$ be another spin structure, with the difference between $s$ and $s_c$ parameterised by $t\in H^1(\mathcal M;\Z_2)$ as in \eqref{eq:spindiff}. According to \eqref{eq:stdresults}, the element $t$ determines a character
\begin{equation} \tilde t\colon \pi_1(\mathcal M)\to \{-1,1\}\subset \U(1),
\end{equation} using the multiplicative action of $\Z_2$ on $\C$. This can be used to construct the line bundle $L_{\tilde t}$ and hence a new spinor bundle
\begin{equation} \mathcal W^t=L_{\tilde t}\otimes \mathcal W^0\end{equation}
corresponding to spin structure $s$. This gives a bundle of spinors $\mathcal W^t$ that is periodic or anti-periodic along a loop according to the value of $\tilde t$. Note that there may not be a spin frame field defined over the whole of $\mathcal M$.

In all these cases the spin frame on the universal covering space can be taken to be $f=1$ and so the formula for the Dirac operator reduces to
\begin{equation}D\psi=\sum_{i=1}^{m+k} \gamma^i \nabla_{v_i}\psi.\label{eq:diracb}
\end{equation}

\begin{example}\label{ex:circle}The circle $\R/2\pi\Z$ has a framing $v_1=\pd{}\theta$. The one-dimensional Clifford algebra has two irreducible representations in $\C$ given by $\gamma^1= i$ or $-i$. The bundle of spinors can be formed from $\R\times\C$ by one of the two quotients $\psi(\theta+2\pi)=\pm\psi(\theta)$, called periodic or anti-periodic boundary conditions. The identification maps  $\pm1\in\Spin(1)$ parameterise the two spin structures; $+1$ is called the Lie group spin structure and $-1$ is called the bounding spin structure. The Dirac operator is $D=\gamma^1\pd{}\theta$.
\end{example}

\begin{example} \label{ex:circle2} Augmenting the circle with one extra dimension leads to the following Dirac operator depending on 
a fixed integer $n$. The orthonormal frame $\R^2\to TS^1_\theta\oplus\R\cong\R^2$ is
\begin{equation}e_\theta=\begin{pmatrix}\cos n\theta &\sin n\theta\\-\sin n\theta &\cos n\theta\end{pmatrix},
\end{equation}
and hence the two vector fields
\begin{equation}v_1=\cos n\theta\, \pd{}\theta,\quad v_2=\sin n\theta\, \pd{}\theta.
\end{equation}
Using the gamma matrices
\begin{equation}\gamma^1=\begin{pmatrix}0&i\\i&0\end{pmatrix},\quad\gamma^2=\begin{pmatrix}0&1\\-1&0\end{pmatrix},\label{eq:2dgamma}\end{equation}
the Dirac operator is
\begin{equation}D= \gamma^1 \nabla_{v_1}+\gamma^2\nabla_{v_2}=\begin{pmatrix}0&ie^{-in\theta}\\ie^{in\theta}&0\end{pmatrix}\nabla_{\pd{}\theta},
\end{equation}
and the chirality operator is 
\begin{equation}\gamma=i\gamma^1\gamma^2=\begin{pmatrix}1&0\\0&-1\end{pmatrix}.\end{equation}
As in Example \ref{ex:circle}, the bundle of spinors can be constructed from $\R\times\C^2$ with either periodic or anti-periodic boundary conditions. 
\end{example}

It is worth noting that the $n=0$ cases of Example \ref{ex:circle2} are the same as the direct sum of the two irreducible cases of Example \ref{ex:circle} with the same boundary conditions, after a change of basis in $\C^2$.

\subsection{Gauge transformations}\label{sec:gauge}

A gauge transformation is a function $R\colon \mathcal M\to\SO(m+k)$ that determines a new frame field $e'=e R$. The gauge transformation determines the homomorphism 
\begin{equation}R_\bullet\colon\pi_1(\mathcal M)\to\pi_1(SO(m+k)).
\end{equation}
If $m+k\ge3$ then $\pi_1(SO(m+k)\cong \Z_2$. In the case when $m+k<3$, then $R_\bullet$ can be composed with the stabilisation $\pi_1(SO(m+k))\to \pi_1(SO(3))\cong\Z_2$ given by the inclusion $SO(m+k)\hookrightarrow SO(3)$. In either case one ends up with a homomorphism
\begin{equation}R'_\bullet\colon\pi_1(\mathcal M)\to\Z_2.
\end{equation}
If $R'_\bullet$ is the trivial homomorphism, $R$ lifts to a function $Z\colon \mathcal M\to\Spin(m+k)$. For a general $R$, there is a lifting of it to a map on the universal covering space $Z\colon \widehat {\mathcal M}\to \Spin(m+k)$.

Let $\mathcal W^0$ be the trivial spinor bundle for frame $e$ using Definition \ref{def:trivial}. Define $\mathcal W'^0$ to be the trivial spinor bundle for frame $e'$ and then define $\mathcal W'^t=L_{\tilde t}\otimes \mathcal W'^0$ with $\tilde t=(-1)^{R'_\bullet}$. Let $\mu$ be the Clifford multiplication of $\mathcal W^0$ and $\mu'$ the Clifford multiplication of $\mathcal W'^t$.

\begin{lemma}\label{lem:Clifford}The map $\Sigma(\mathcal W'^t)\to\Sigma(\mathcal W^0)$ defined by $\psi\mapsto Z\psi$ intertwines the Clifford multiplications, i.e.,
\begin{equation} \mu'=Z^{-1}\,\mu\,\left(1\otimes Z\right).
\end{equation}
\end{lemma}

\begin{proof} The computation is done on the covering space, where both bundles are trivial and the spin frames are $f=1$ and $f'=1$.
In the following equations, the point $x$ is omitted to simplify the notation. 
\begin{multline}\mu'(\nu\otimes w)=\rho'(\nu)w=\check\rho'(P^*\nu)w\\
=\rho_0(e'^*P^*\nu)w=\rho_0(R^*e^*P^*\nu)w=Z^{-1}\rho_0(e^*P^*\nu)Zw\\
=Z^{-1}\check\rho(P^*\nu)Zw
=Z^{-1}\rho(\nu)Zw=Z^{-1}\mu(\nu\otimes Zw).\end{multline}
Note that the third equality uses \eqref{eq:gencliff} and the fifth equality uses \eqref{eq:intertwineZ}.
\end{proof}
The two spinor bundles have the same spin structure. The line bundle in the construction of $\mathcal W'^t$ cancels the change in homotopy class of the orthonormal frame field.

\begin{example} In Example \ref{ex:circle2}, a gauge transformation is defined by 
\begin{equation} R(\theta)=\begin{pmatrix}\cos j\theta &\sin j\theta\\-\sin j\theta &\cos j\theta\end{pmatrix}
\end{equation}
for a fixed integer $j$. The two  lifts to $\Spin(2)\cong\U(1)$ are $\pm e^{i j\theta/2}$, represented in the Clifford module $\C^2$ by
\begin{equation}Z(\theta)= \pm\begin{pmatrix} e^{i j\theta/2}&0 \\ 0&e^{-i j\theta/2}.\end{pmatrix}
\end{equation}
If $l$ denotes the circle, then $R'_\bullet(l)=j \pmod 2$. If $j$ is odd, $Z$ maps periodic spinors into anti-periodic spinors, and vice-versa.

Applying the gauge transformation to the case $n$ of Example \ref{ex:circle2} takes it to the case $n+j$, and in particular, it relates the case $n=0$ to the case $n=j$. This shows that the spin structure is the Lie group structure if $n$ is even and the boundary conditions are periodic, or $n$ is odd and the boundary conditions anti-periodic. Otherwise it is the bounding spin structure.
\end{example}

In general, the Dirac operator is preserved under a gauge transformation. More precisely, this is the following result.

\begin{lemma}\label{lem:gauge} If $D$ is the Dirac operator on $\Sigma(\mathcal W^0)$ and $D'$ the Dirac operator on $\Sigma(\mathcal W'^t)$, then  $D'=Z^{-1}DZ$.
\end{lemma}

\begin{proof} A connection on the orthonormal frame bundle $\mathcal O$ determines a connection on any spin frame bundle uniquely, and hence a covariant derivative on the associated bundle of spinors. The two spin frame bundles used to construct $\mathcal W^0$ and $\mathcal W'^t$ are isomorphic and therefore so are the associated spinor bundles, with the isomorphism given in Lemma \ref{lem:Clifford}. Hence the covariant derivative on a section of $\mathcal W'^t$ in the direction of vector $\xi$  is given by $\nabla'_\xi\psi'=Z^{-1}\nabla_\xi Z\psi'$, or more abstractly
\begin{equation}\nabla'=(1\otimes Z^{-1})\nabla Z.
\end{equation}
Then, using Lemma \ref{lem:Clifford},
 \begin{equation}D'=\mu'\nabla'=Z^{-1}\mu(1\otimes Z)\;
(1\otimes Z^{-1})\nabla Z=Z^{-1}DZ.
\end{equation}
\end{proof}


\subsection{Dirac Operator on the Torus}\label{sec:diractorus}

The description of the Dirac operator on a manifold is specialised to the case of a torus. The spectrum of the Dirac operator is calculated for the general flat torus using a parallel frame.

Let $A\in\GL(2,\R)$. Then a general parallel frame field on the torus is determined by the constant matrix 
\begin{equation}A^{-1}=\frac1{ad-bc}\begin{pmatrix}d&-b\\-c&a\end{pmatrix}.
\end{equation}
The spin structure is determined by $\sigma\in\Z_2\times\Z_2\cong H^1(T^2;\Z_2)$, as in Example \ref{ex:torusss}. The spinor bundle is $\mathcal W^\sigma$, as constructed in Section \ref{sec:framefields}. This means that $\sigma$ determines whether the spinor fields are periodic or anti-periodic along each of the two axes of the torus. 

The two cases of spinors of interest here are the irreducible case $p=2$ and the case $p=4$, where $ p $ is the dimension of the Clifford module.
The Dirac operator is 
\begin{equation}D= \frac1{ad-bc}\left( \gamma^4( d\partial_{\theta}-c\partial_{\phi}) + \gamma^2(-b \partial_\theta+a\partial_\phi)\right),
\end{equation}
using two gamma matrices (labelled with $2$ and $4$ for later convenience).
The square of $D$ is the Laplacian \eqref{eq:laplacian}. The spectrum is calculated by substituting the plane waves
\begin{equation}
\psi(\theta,\phi) = 
 e^{i(k \theta + l \phi)}\,\psi_0\label{eq:planewaves}
\end{equation}
into the Dirac equation and solving for the spinor $\psi_0\in\C^p$. The possible values of $k$ and $l$ are either integer or half-integer and are determined by the spin structure, 
\begin{equation}2(k,l)=\sigma \pmod 2.\label{eq:sskl}\end{equation}
The eigenvalue is given by the length of the vector $(k,l)$ using the appropriate metric \cite{Bar},
\begin{equation}\lambda_{k,l,\pm}=\pm\frac1{ad-bc}\sqrt{(dk-cl)^2+(al-bk)^2},\label{eq:commtorusev}
\end{equation}
for $(k,l)\ne(0,0)$. 

In the irreducible case with $p=2$, this eigenvalue has multiplicity one for given values of $k$, $l$ and $\pm$.  
For the spin structure $(0,0)$, there is also the eigenvalue $\lambda_{0,0}=0$ with multiplicity two. 

The case $p=4$ is the one of interest below. This occurs by introducing two additional gamma matrices $\gamma^1$ and $\gamma^3$, as used for dimension $m=4$. This can be thought of as arising from the tensor product of the two-dimensional spinors with another trivial two dimensional spinor module, $\C^4\cong\C^2\otimes\C^2$. For this case $p=4$, all the above multiplicities are doubled.

The spinors have a real structure $j\colon \C^p\to\C^p$ that commutes with the gamma matrices and satisfies the conditions \eqref{eq:table}. Extending $j$ to an antilinear map $J$ on spinor fields, it commutes with the Dirac operator, $JD=DJ$. If $\psi$ is an eigenvector, then so is $J\psi$, which has the same eigenvalue. This can be seen from  \eqref{eq:commtorusev}, since $\lambda_{k,l,\pm}=\lambda_{-k,-l,\pm}$.

The torus admits a $ U(1) \times U(1) $ action by isometries, translating along each axis. This action is covered by an action of $\Spin(2)\times\Spin(2)$ on the spinor bundle. The elements of this group are parameterised by $(\Theta, \Phi)$, each variable having period $4\pi$. The action on the spinor fields is
\begin{equation}
\Pi_{( \Theta, \Phi)}\, \psi ( \theta, \phi) = \psi(\theta + \Theta, \phi + \Phi).
\end{equation} 
This commutes with the Dirac operator and the plane waves in \eqref{eq:planewaves} are eigenvectors of $P_{ (\Theta, \Phi)}$ with eigenvalue $e^{i(k\Theta+l\Phi)}$.

\subsection{Square torus}\label{sec:squaretorus}
 The rotating frame field that generalises to the fuzzy case (in Section \ref{sec:fuzzydirac}) is introduced, giving a formula for the Dirac operator using this frame field. This is done first for the unit square torus, which is the case $a=d=1$, $b=c=0$. 
The rotating frame field uses the tangent space augmented with $\R^2$. Thus there are four gamma matrices acting on spinors in $\C^4$, and the real structure $j$ obeys the relations for dimension four in \eqref{eq:table}. 

The Dirac operator with the parallel frame specialises to
\begin{equation}D_1=\gamma^2\partial_\phi+\gamma^4\partial_\theta.
\end{equation}
The parallel frame  $e\colon\R^4\to T(T^2)\oplus\R^2$ is defined by $e(\xi_2)=\partial_\phi$, $e(\xi_4)=\partial_\theta$, and $e(\xi_1)=(1,0)$, $e(\xi_3)=(0,1)$ are basis vectors in the additional subspace $\R^2$. Applying the $\SO(4)$ gauge transformation
\begin{equation}
R(\theta,\phi) = \left( \begin{array}{cccc}
\cos\theta & \sin\theta & 0 & 0 
\\
-\sin\theta & \cos\theta & 0 & 0 
\\
0 & 0 & \cos\phi & \sin\phi
\\
0 & 0 & -\sin\phi & \cos\phi
\end{array} \right)
\end{equation}
gives a new frame $e'=eR$.
The vector fields on $T^2$ for the frame $e'$ are then
\begin{equation} \begin{aligned}\label{eq:rotatingvectors}
w_1 = -(\sin\theta) \partial_{\phi} &\quad w_2 = (\cos\theta) \partial_{\phi}
\\
w_3 = -(\sin\phi) \partial_{\theta} &\quad w_4 = (\cos\phi) \partial_{\theta}
\end{aligned}
\end{equation}
A lift of $R$ to $\Spin(4)$ is given by
\begin{equation} \begin{aligned}
Z_1 = \exp\left(-\frac{1}{2}\left( \theta\gamma^1 \gamma^2 +  \phi\gamma^3 \gamma^4\right)\right).
\end{aligned}
\end{equation}
Due to the factor of $1/2$ in this formula, $Z$ is anti-periodic in both $\theta$ and $\phi$. Therefore it determines a map 
of $\mathcal W^\sigma$ to the bundle $\mathcal W'^{\sigma+(1,1)}$, following the discussion in Section \ref{sec:gauge}.

Using Lemma \ref{lem:gauge}, the Dirac operator on $T^2$ with the rotating frame is
\begin{equation} \begin{aligned}
 D'_1=Z_1^{-1}D_1Z_1 =& \bigl(-(\sin\theta) \gamma^1 + (\cos\theta)\gamma^2\bigr) \bigl( \partial_{\phi}  -\frac{1}{2} \gamma^3 \gamma^4 \bigr) 
\\ &+ \bigl(-(\sin\phi)  \gamma^3 + (\cos\phi)\gamma^4 \bigr) \bigl(\partial_{\theta} -\frac{1}{2} \gamma^1 \gamma^2 \bigr). \end{aligned}\label{eq:rotatingD}
\end{equation}
Explicit formulas are given by choosing the gamma matrices 
\begin{equation}
\gamma^1 = \left( \begin{array}{cccc}
0 & 0 & i & 0
\\
0 & 0 & 0 & i
\\
i & 0 & 0 & 0
\\
0 & i & 0 & 0
\end{array} \right) \quad
\gamma^2 = \left( \begin{array}{cccc}
0 & 0 & -1 & 0
\\
0 & 0 & 0 & 1
\\
1 & 0 & 0 & 0
\\
0 & -1 & 0 & 0
\end{array} \right) 
$$$$
\gamma^3 = \left( \begin{array}{cccc}
0 & 0 & 0 & i
\\
0 & 0 & -i & 0 
\\
0 & -i & 0  & 0 
\\
i & 0 & 0 & 0
\end{array} \right) 
\quad  \gamma^4 = \left( \begin{array}{cccc}
0 & 0 & 0 & -1
\\ 
0 & 0 & -1 & 0
\\
0 & 1 & 0 & 0
\\ 
1 & 0 & 0 & 0
\end{array} \right).\label{eq:gamma}
\end{equation}
The operator $D'_1$ has eigenvectors
\begin{equation}
\psi_{k',l',\pm}^1 = 
\left( \begin{array}{c} 
e^{i(k'+1) \theta + i (l'+1)\phi}
\\
e^{i k' \theta +il'\phi}
\\
\pm i\frac{l'+k'+1}{\sqrt{(k'+1/2)^2+(l'+1/2)^2}} e^{i k' \theta + i (l'+1) \phi}
\\
\pm i \frac{(k'-l')}{\sqrt{(k'+1/2)^2+(l'+1/2)^2}}  e^{i (k'+1)\theta + il'\phi}
\end{array} \right)
\end{equation}
and
\begin{equation}
\psi_{k',l',\pm}^2 = 
\left( \begin{array}{c} 
\pm i\frac{l'+k'+1}{\sqrt{(k'+1/2)^2+(l'+1/2)^2}} e^{i(k'+1) \theta + i (l'+1)\phi}
\\
\pm i \frac{(k'-l')}{\sqrt{(k'+1/2)^2+(l'+1/2)^2}} e^{i k' \theta +il'\phi}
\\
e^{i k' \theta + i (l'+1) \phi}
\\
e^{i (k'+1)\theta + il'\phi}
\end{array} \right),
\end{equation}
with $ (2k'+1,2l'+1) = \sigma \pmod 2$. The corresponding eigenvalues are given by
\begin{equation}
 \pm \sqrt{(k'+1/2)^2 + (l'+1/2)^2},  \label{eq:classevBB}
\end{equation}
each with multiplicity 2. This agrees with \eqref{eq:commtorusev} by setting $k'+1/2=k$ and $l'+1/2=l$.  Calculating the square of the rotating Dirac operator gives
\begin{equation}\begin{aligned}
{D_1'}^2 &= Z_1^{-1}D^2_1Z_1 = Z_1^{-1}(-\partial_\phi^2-\partial_\theta^2)Z_1\\
&=-\partial_{\theta}^2 - \partial_{\phi}^2 + \gamma^1 \gamma^2 \partial_{\theta} + \gamma^3 \gamma^4 \partial_{\phi} + \frac{1}{2}. 
\end{aligned}\end{equation}
The covariant derivative on spinors determines a Laplace operator on spinors  \cite[Theorem 8.8]{Lawson}. In general, the Lichnerowicz-Schr\"{o}dinger equation relates the square of the Dirac operator to this Laplacian and the curvature scalar. The curvature scalar is zero for a flat torus, so in this case the square of the Dirac operator is equal to the Laplace operator on spinors.


\subsection{Torus with integral metric}\label{sec:commintegral} The case of the Dirac operator on a torus with a general integral metric and rotating frame is now discussed. Consider the transformation $ A\colon( \theta, \phi ) \mapsto  (a \theta +b \phi, c\theta+d\phi) $. This transformation is used to pull back functions on the square torus, giving
\begin{equation}
C \mapsto U = e^{i(a \theta + b \phi)},\quad \quad S \mapsto V = e^{i(c \theta + d \phi)}.
\end{equation}
Vector fields are also related by pull-back, so for vector fields of the rotating frame \eqref{eq:rotatingvectors},  $ v'_i=A^\bullet w_i=A^{-1}_\bullet w_i$. Explicitly, these are
\begin{equation}
\begin{aligned}
&v'_1 = \frac{1}{a d-b c}(b \sin(a \theta + b \phi) \partial_{\theta} - a \sin(a \theta + b \phi) \partial_{\phi})
\\
&v'_2 = \frac{1}{a d-b c}(-b \cos(a \theta + b \phi) \partial_{\theta} + a \cos(a \theta + b \phi) \partial_{\phi})
\\
&v'_3 = \frac{1}{a d-b c}(-d \sin(c \theta + d \phi) \partial_{\theta} + c \sin(c \theta + d \phi) \partial_{\phi})
\\
&v'_4 = \frac{1}{a d-b c}(d \cos(c \theta + d \phi) \partial_{\theta} -c \cos(c \theta + d \phi) \partial_{\phi}).
\end{aligned}
\end{equation}
Therefore, the lift to the action induced on the spin bundle $\mathcal W^\sigma$ is
\begin{equation}\label{eq:spinlift}
\begin{aligned}
Z = \exp\left(-\frac{1}{2}\left(\gamma^1 \gamma^2 (a\theta+ b \phi) + \gamma^3 \gamma^4(c \theta +d \phi)\right)\right ),
\end{aligned}
\end{equation}
so that $\psi'=Z^{-1}\psi$. The Dirac operator transforms to
\begin{equation}
\begin{aligned}
D' &= Z^{-1}D Z\\
&= \frac{1}{ad-bc}\bigl(b(\sin(a \theta + b \phi)\gamma^1 - \cos(a \theta + b \phi)\gamma^2) -d( \sin(c \theta + d \phi)\gamma^3 
\\&- \cos(c \theta + d \phi) \gamma^4) \bigr) \bigl( \partial_{\theta} - \frac{1}{2}(a \gamma^1 \gamma^2 + c \gamma^3 \gamma^4)\bigr)
\\
&+\frac{1}{ad-bc}\bigl(-a(\sin(a \theta + b \phi) \gamma^1 - \cos(a \theta + b \phi) \gamma^2) + c( \sin(c \theta + d \phi) \gamma^3
\\
&-\cos(c \theta + d \phi) \gamma^4)\bigr)\bigl(\partial_{\phi} - \frac{1}{2}(b \gamma^1 \gamma^2 + d \gamma^3 \gamma^4)\bigr)
\end{aligned}
\end{equation}
acting on sections of $\mathcal W'^{\sigma'}$, with $\sigma'=\sigma+\bigl(\lBracket a\rBracket_2+\lBracket c\rBracket_2,\lBracket b\rBracket_2+\lBracket d\rBracket_2\bigr)$.
It is easy to check that one recovers $D_1^{\prime} $ when $a=d=1,b=c=0$. 

The transformation $A$ acts on the spinor fields by regarding them as $\C^4$-valued functions on $\R^2$ and using the pull-back of functions, denoted $A^\bullet$. (The spinor components are not transformed.) The Dirac operator $D'$ can be characterised by this pull-back. It obeys
\begin{equation} D' A^\bullet=A^\bullet D_1'.\end{equation}

The spinor Laplacian for the general integral torus in the rotating frame is 
\begin{equation}
\begin{aligned}
(D')^2 &= \frac{-1}{(ad-bc)^2}\left((b^2+d^2)\pd{\mathstrut^2}{\theta^2}-2(ab+cd)\pd{\mathstrut^2}{\theta\partial\phi}+(a^2+c^2)\pd{\mathstrut^2}{\phi^2}\right)
\\ &+ \frac{1}{ad-bc} (\gamma^1 \gamma^2 (d \partial_{\phi} - c  \partial_{\theta} ) + \gamma^3 \gamma^4 (a \partial_{\phi}- b \partial_{\theta} ) ) + \frac{1}{2}.
\end{aligned} \label{eq:diracsquare}
\end{equation}

The $\Spin(2)\times \Spin(2)$ action in the rotating frame becomes
\begin{equation}
\Pi^{\prime}_{(\Theta,\Phi)} \psi'(\theta,\phi) = W \psi'(\theta + \Theta,\phi + \Phi),  \label{eq:pispinaction}
\end{equation} 
with a non-trivial gauge transformation  
\begin{equation}\label{eq:Wspinaction}
W=Z(\Theta,\Phi) = \text{exp} \left( -\frac{1}{2} (\gamma^1 \gamma^2 ( a \Theta + b \Phi) + \gamma^3 \gamma^4 ( c \Theta + d \Phi )  ) \right)
\end{equation}
on the spinors. The Dirac operator $ D^{\prime} $ then commutes with the action of $\Spin(2)\times \Spin(2)$ on the spinor fields, i.e.,
\begin{equation}
D^{\prime} \Pi^{\prime}_{(\Theta,\Phi)}  = \Pi^{\prime}_{(\Theta,\Phi)} D^{\prime}.
\end{equation}

\section{Dirac operator on the fuzzy torus}\label{sec:fuzzydirac}
 
The Dirac operator on a fuzzy space is the fundamental structure that encodes the geometry of the space. It is characterised by a set of algebraic axioms for a structure called a finite real spectral triple.  This section introduces a particular spectral triple for each fuzzy torus. This gives a Dirac operator for the non-commutative analogues of the tori with integral metrics and arbitrary spin structures, as presented in Section \ref{sec:commintegral}.

Real spectral triples were originally defined in \cite{Connes:1995tu} and \cite{Connes:1996gi}, the axioms undergoing some slight modifications since then. The definition includes compact spin manifolds in the commutative case. Specialising to the finite case simplifies the definition because the analytic axioms are not required. The precise definition of the finite case is the one given in  \cite{Barrett:2015naa} and is summarised briefly here for the even-dimensional cases.

\begin{definition}
A finite real spectral triple is
\begin{itemize}
\item An integer $d$ defined mod 8, called the KO-dimension.
\item A finite-dimensional Hilbert space $\mathcal H$.
\item A $\ast$-algebra $\mathcal A$ of operators with a faithful left action in $\mathcal H$.
\item A Hermitian operator $\Gamma\colon  \mathcal H\to\mathcal H$, called the chirality, that commutes with the algebra action and obeys $\Gamma^2=1$.
\item An antiunitary map $\mathcal J\colon \mathcal H\to\mathcal H$, called the real structure, such that
\begin{equation}[\mathcal Ja\mathcal J^{-1},b]=0\end{equation} 
for all $a, b\in\mathcal A$. The real structure satisfies $\mathcal J^2=\epsilon$, $\mathcal J\Gamma \mathcal J^{-1}\Gamma=\epsilon''$, using (for even $d$) the sign table in  Section \ref{sec:diracmfd}.
\item A Hermitian operator $D\colon \mathcal H\to\mathcal H$, called the Dirac operator,  satisfying
\begin{equation}[[D,a],\mathcal Jb\mathcal J^{-1}]=0\label{eq:foc}\end{equation} 
for all $a,b\in\mathcal A$, and (for even $d$), $D\Gamma+\Gamma D=0$ and $D\mathcal J=\mathcal JD$.
\end{itemize}
\end{definition}
\bigskip

The axioms determine a right action of $a\in\mathcal A$ on $\psi\in \mathcal H$ by
\begin{equation} \psi\, a= \mathcal Ja^*\mathcal J^{-1}\psi.\end{equation}
This commutes with the left action and so makes $\mathcal H$ a bimodule.

 The fuzzy analogue presented here is of a flat torus with a tangent space augmented with two extra dimensions and a non-constant framing. 
This uses the space of spinors $\C^4$ regarded as a Clifford module of type $(0,4)$. These spinors are twice the dimension of the usual Dirac spinors on a torus. The chirality operator for them is $\gamma = \gamma^1 \gamma^2 \gamma^3 \gamma^4$ and the real structure is
\begin{equation*}
j \left( \begin{array}{c}
v_1
\\
v_2
\\
v_3
\\
v_4
\end{array} \right) = \left( \begin{array}{c}
\overline{v_2}
\\
-\overline{v_1}
\\
-\overline{v_4}
\\
\overline{v_3}
\end{array} \right),
\end{equation*}
which satisfies $j\gamma^i=\gamma^i j$.

The following definition is phrased in a general way but all of the examples that follow will have a very simple construction in terms of spaces of matrices.

\begin{definition}\label{def:spectral}
Let $U,V,\mathfrak h$ be a fuzzy torus with real structure $J\colon\mathfrak h\to\mathfrak h$. Also, let $ X, Y \in \langle U, V \rangle $ with $ XY = \X YX $ and choose a fourth root $\X^{1/4}$. The spectral triple for the fuzzy torus is as follows.
\begin{itemize}
\item The KO-dimension is 4.
\item The Hilbert space is $\mathcal H=\C^4\otimes \mathfrak h$.
\item The $\ast$-algebra is $\mathcal A=\langle U,V\rangle$, acting in $\mathcal H$ on the left by
\begin{equation*} U(v\otimes m)=v\otimes Um ,\quad \quad  V(v\otimes m)=v\otimes Vm
\end{equation*}
 \item The real structure is $\mathcal J(v\otimes m)=jv\otimes J m$. 
\item The chirality operator is $\Gamma(v\otimes m)=\gamma v\otimes m$.
\item The Dirac operator on the fuzzy torus is
\begin{equation*} \begin{aligned}
&D_{X,Y} = -\frac{1}{4(\X^{1/4} - \X^{-1/4})} (\gamma^1 \otimes [ X + X^{*},{} \cdot{} ]+ \gamma^2 \otimes  i[X^{*} - X , {}\cdot{} ] 
\\ &-  \gamma^3 \otimes [Y+Y^{*},{}\cdot{}] - \gamma^4 \otimes i[Y^{*} - Y,{} \cdot{} ] )
\\&- \frac{1}{4(\X^{1/4} +\X^{-1/4})}  ( \gamma^2 \gamma^3 \gamma^4 \otimes \{X+X^{*},{}\cdot{}\} - \gamma^1 \gamma^3 \gamma^4 \otimes i \{X^{*} - X,{}\cdot{} \} 
\\&+ \gamma^1 \gamma^2 \gamma^4 \otimes \{ Y + Y^{*}, {}\cdot{} \} - \gamma^1 \gamma^2 \gamma^3 \otimes i \{Y^{*} - Y, {}\cdot{} \}),
\end{aligned}\end{equation*}
using the right action determined by $J$.
\end{itemize} 
\end{definition}

It is a straightforward calculation to check that this defines a real spectral triple. In particular, the first order condition \eqref{eq:foc} follows from the fact that $D$ is a sum of terms that commute with either the left action or the right action of $U$ and $V$. 

The Dirac operator in Definition \ref{def:spectral} may be written a little more systematically, as in  \cite{Barrett:2015naa}, as
\begin{equation}
\begin{aligned}
D_{X,Y} = \frac{1}{\X^{1/4}-\X^{-1/4}}\sum_{i} \gamma^i \otimes [K_i, \cdot\, ] + \frac{1}{\X^{1/4}+\X^{-1/4}}\sum_{i<j<k} \gamma^i \gamma^j \gamma^k \otimes \{K_{ijk},\cdot\}
\end{aligned}
\end{equation}
with
\begin{equation}
\begin{aligned}
K_{1} = K_{234} = -\frac{1}{4}(X+ X^{*}), & \quad\quad K_2 = -K_{134} = -\frac{i}{4}(X^*-X),
\\
K_{3} = -K_{124} = \frac{1}{4}(Y+Y^*), &\quad\quad K_4 = K_{123} = \frac{i}{4}(Y^*-Y).
\end{aligned}
\end{equation}

 The most delicate aspects of this definition are the numerical coefficients involving $\X^{1/4}$. Some insight into these is obtained by calculating the square of the Dirac operator. A lengthy calculation shows
\begin{equation} \begin{aligned}
&D_{X,Y}^2 = -\frac{1}{4(\X^{1/4}-\X^{-1/4})^2} ( 1 \otimes [ X, [X^{*} , \cdot\, ] + 1 \otimes [Y, [Y^{*},\cdot\,] ) 
\\&+ \frac{i}{2(\X^{1/2}-\X^{-1/2})}(\gamma^1 \gamma^2 \otimes \{Y ,[ Y^{*} ,\cdot\,] \} 
- \gamma^3 \gamma^4 \otimes \{ X, [ X^*, \cdot\, ] \})
\\& +\frac{1}{4(\X^{1/4}+\X^{-1/4})^2}( 1 \otimes \{X,\{X^{*},\cdot\}\} + 1 \otimes \{Y,\{Y^{*},\cdot\} \}).
\end{aligned}\end{equation}

The calculation of the square can be broken down in the following way. Setting
\begin{equation*} \begin{aligned}
&E_{X} = -\frac{1}{4(\X^{1/4} - \X^{-1/4})} (\gamma^1 \otimes [ X + X^{*},{} \cdot{} ]+ \gamma^2 \otimes  i[X^{*} - X , {}\cdot{} ]  )
\\&- \frac{1}{4(\X^{1/4} +\X^{-1/4})}  ( \gamma^2 \gamma^3 \gamma^4 \otimes \{X+X^{*},{}\cdot{}\} - \gamma^1 \gamma^3 \gamma^4 \otimes i \{X^{*} - X,{}\cdot{} \} )
\end{aligned}\end{equation*}
and 
\begin{equation*} \begin{aligned}
&E_{Y} = -\frac{1}{4(\X^{1/4} - \X^{-1/4})} (-  \gamma^3 \otimes [Y+Y^{*},{}\cdot{}] - \gamma^4 \otimes i[Y^{*} - Y,{} \cdot{} ] )
\\&- \frac{1}{4(\X^{1/4} +\X^{-1/4})}  (  \gamma^1 \gamma^2 \gamma^4 \otimes \{ Y + Y^{*}, {}\cdot{} \} - \gamma^1 \gamma^2 \gamma^3 \otimes i \{Y^{*} - Y, {}\cdot{} \}).
\end{aligned}\end{equation*}
Then $D_{X,Y}=E_X+E_Y$. The crucial property is
\begin{lemma}
\begin{equation}E_XE_Y+E_YE_X=0\label{eq:Esanticommute}
\end{equation}
\end{lemma}
\begin{proof}Put $\alpha=-\frac{1}{4(\X^{1/4} - \X^{-1/4})}$, $\beta=- \frac{1}{4(\X^{1/4} +\X^{-1/4})} $ and use the relations 
\begin{equation}\begin{aligned}
[[A,\cdot\,],[B,\cdot\,]]&=[[A,B],\cdot\,]\\
\{[A,\cdot\,],\{B,\cdot\,\}\}+\{\{A,\cdot\,\},[B,\cdot\,]\}&=2[\{A,B\},\cdot\,]\\
[\{A,\cdot\,\},\{B,\cdot\,\}]&=[[A,B],\cdot\,]
\end{aligned}\end{equation}

The coefficient of $\gamma^1\gamma^3$ is
\begin{multline}-\alpha^2[[X+X^*,Y+Y^*],\cdot\,]
-\beta^2[[X+X^*,Y+Y^*],\cdot\,]\\
-2\alpha\beta [\{i(X^*-X),i(Y^*-Y)\},\cdot\,]
\end{multline}
This is the commutator with the negative of
\begin{multline}(\alpha^2+\beta^2)[X+X^*,Y+Y^*]+2\alpha\beta \{i(X^*-X),i(Y^*-Y)\}\\
=(\alpha^2+\beta^2)((Q-1)(YX+Y^*X^*)+(Q^*-1)(YX^*+Y^*X))\\
-2\alpha\beta(Q+1)(YX+Y^*X^*-(Q^*+1)(YX^*+Y^*X))
\end{multline}
which is zero due to the relation
\begin{equation}
(\alpha^2+\beta^2)(Q-1)-2\alpha\beta(Q+1)=0.\end{equation}
The coefficients of other pairs of gamma matrices vanish by a similar calculation.
\end{proof}
The lemma implies that $D_{X,Y}^2=E_X^2+E_Y^2$.
The property \eqref{eq:Esanticommute} does actually fix the numerical coefficients uniquely, up to an overall constant.  It is worth noting that the formulas for $E_X$ and $E_Y$ correspond if one interchanges $\gamma^1\leftrightarrow\gamma^3$, 
$\gamma^2\leftrightarrow\gamma^4$, $X\leftrightarrow Y$ and $\X^{1/4}\leftrightarrow \X^{-1/4}$.

\subsection{Square Fuzzy Torus}
To calculate explicit features of the fuzzy torus Dirac operator one must consider specific choices of $X$ and $Y$. This corresponds to picking a specific metric on a fuzzy torus. Mirroring Section \ref{sec:dirac}, the square fuzzy torus is considered first. This is the simplest case, allowing a discussion of the main features of the Dirac operator without introducing a lot of formalism. More general cases are discussed in Sections \ref{sec:fuzzyint} and \ref{sec:fuzzyspin}.

Let  $ \mathfrak h = M_N(\mathbb{C}) $, $J = \ast $, $X=U=C$ and $Y=V=S$, so that  $\mathcal A= M_N(\C)$.

The first thing to note is that the Dirac operator corresponds to the commutative case \eqref{eq:rotatingD} using the heuristics \eqref{eq:commutatorlimit} for the commutators, adding
\begin{equation}
\begin{gathered}\{C,{\cdot}\}\rightsquigarrow 2e^{i\theta},\quad \{C^*,{\cdot}\}\rightsquigarrow 2e^{-i\theta}\\
\{S,{\cdot}\}\rightsquigarrow 2e^{i\phi},\quad \{S^*,{\cdot}\}\rightsquigarrow 2e^{-i\phi}\end{gathered} \label{eq:anticomm}
\end{equation}
for the anticommutators, then setting $\q^{1/4}$ to $1$. 

Define the vector
\begin{equation}
\psi = \left( \begin{array}{c}
\alpha_{kl}  \e^{(k+1,l+1)}
\\
\beta_{kl}  \e^{(k,l)} 
\\
\gamma_{kl} \e^{(k,l+1)} 
\\
\delta_{kl} \e^{(k+1,l)} 
\end{array} \right),
\end{equation}
with $\alpha_{kl},\beta_{kl},\gamma_{kl},\delta_{kl} $ fixed complex numbers, each a function of integers $k$ and $l$. 
This is an eigenvector of the Dirac operator, the eigenvalue equation being
 \begin{equation}
D_{C,S}^2 \psi = ([k+1/2]_\q^2 + [l+1/2]_\q^2)\psi.\label{eq:CSev}
\end{equation}
Each eigenvalue has multiplicity 4, as expected from the spinor doubling. 
It follows immediately that the eigenvalues of $D_{C,S}$ are given by
\begin{equation}
\lambda_{k,l,\pm} = \pm \sqrt{[k+1/2]_\q^2 + [l+1/2]_\q^2} \label{eq:fuzzyevBB}
\end{equation}
with each eigenvalue having multiplicity 2. 

\begin{figure}
\begin{center} 
\includegraphics[scale=0.55]{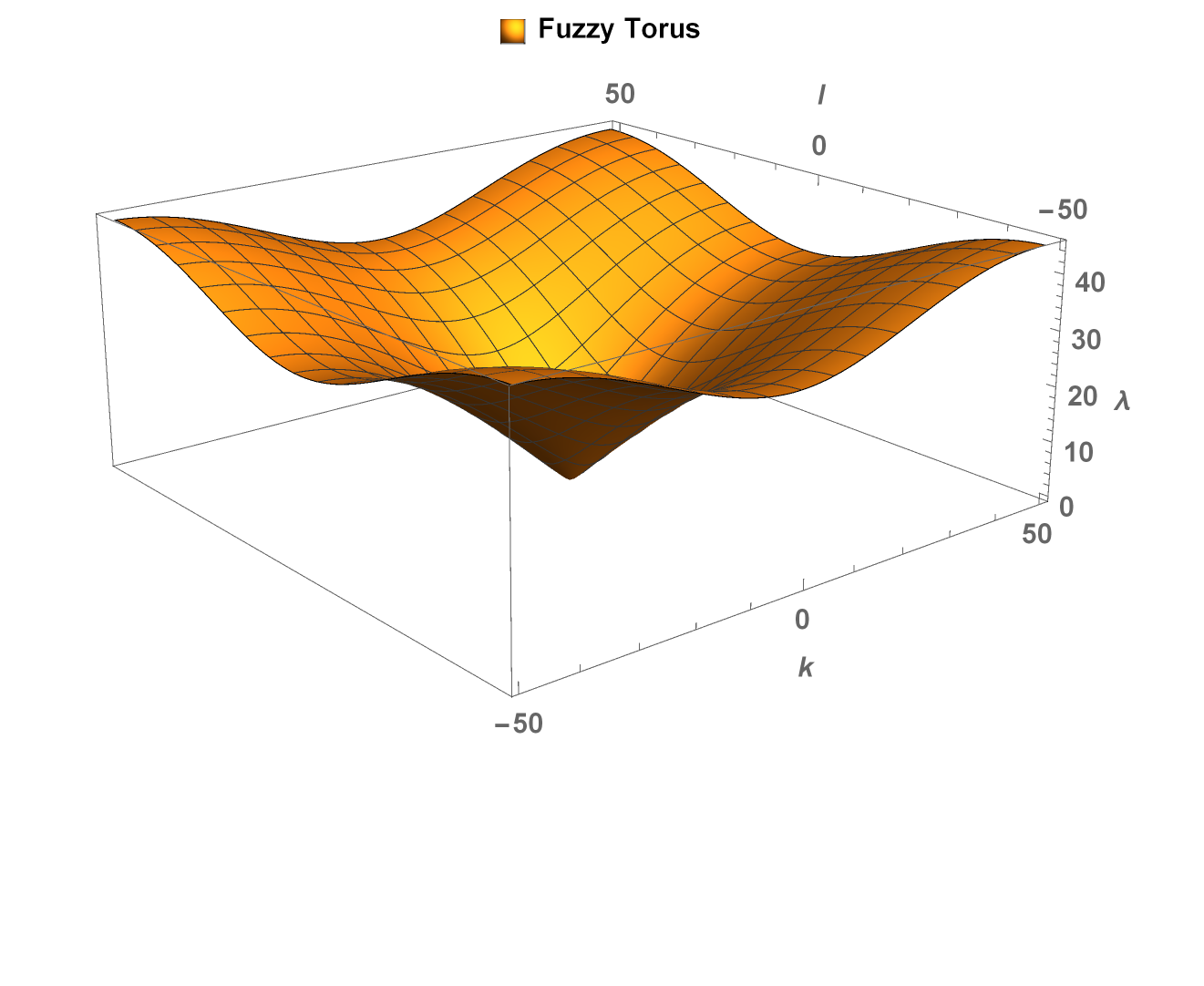}
\end{center}
\vspace{-2.5cm}
\caption{A 3D plot of the positive eigenvalues \eqref{eq:fuzzyevBB} with $N=100$. }\label{fig:3dplot}
\end{figure}

The eigenspaces are periodic in the labels $k$, $l$, and so a unique labelling is given by
\begin{equation}-\frac N2<k,l\le \frac N2.\label{eq:cutoff}
\end{equation}
Plotting \eqref{eq:fuzzyevBB} for a fixed value of $N$ shows the geometry of the fuzzy torus. As in Section \ref{sec:fuzzy},  the value $\q = e^{2 \pi i/N}$ is used for all of the following plots. Figure \ref{fig:3dplot} shows the dispersion relation, the discrete values for integer $k$ and $l$ being interpolated by a smooth surface.  This surface
 would be periodic if $k$ and $l$ were continued to $\Z\times \Z$, and looks like a dispersion relation for a massless electron in a lattice in solid state physics. 

\begin{figure}
\begin{center} 
\includegraphics[scale=0.5]{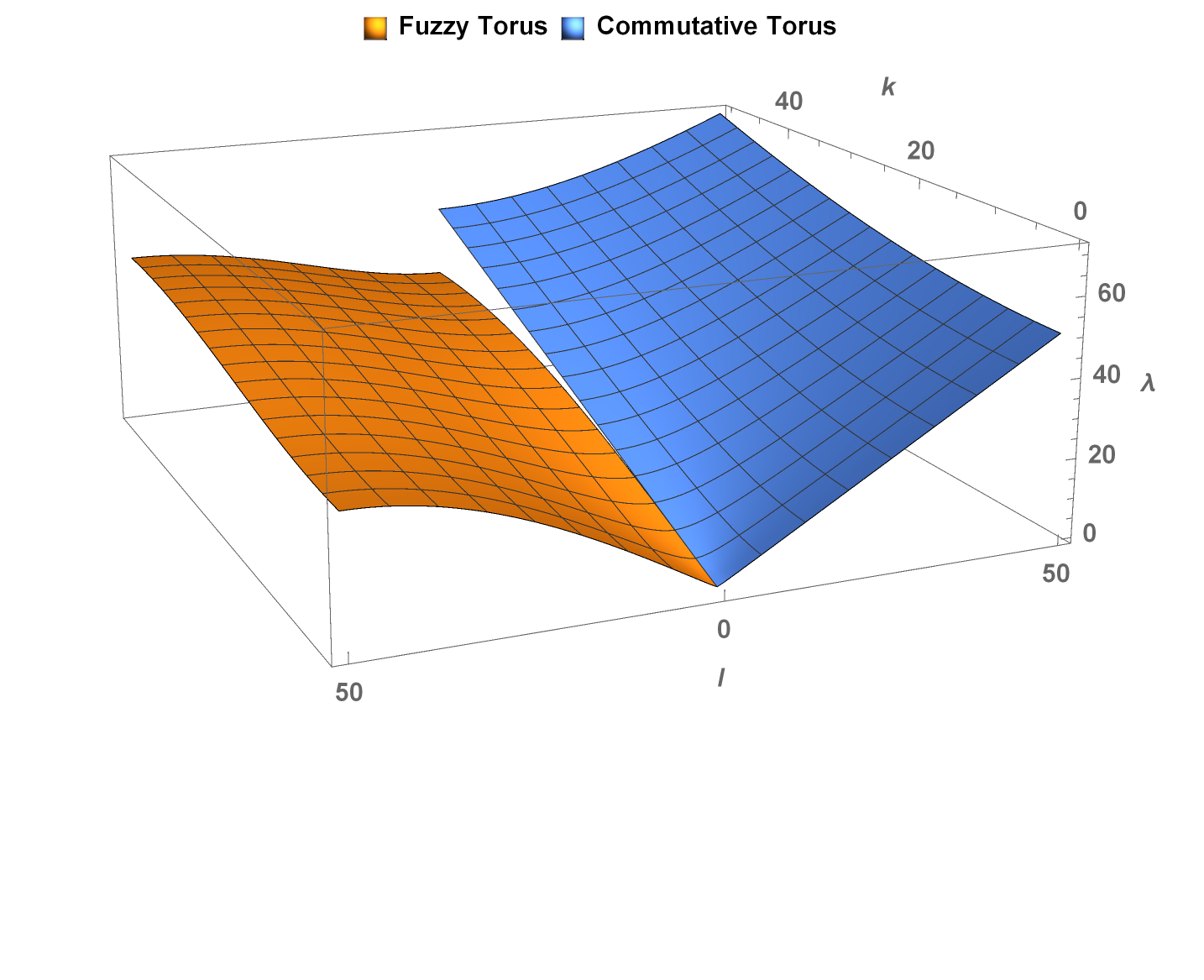}
\end{center}
\vspace{-2.5cm}
\caption{A side by side comparison of the positive eigenvalues \eqref{eq:fuzzyevBB} and \eqref{eq:classevBB}, with $N=100$ and $k,l>0$.  }\label{fig:3dcomp}
\end{figure}

Taking the limit $\q^{1/2}\rightarrow 1 $ with $k $ and $ l $ fixed gives
\begin{equation}
[k+1/2]_\q^2 + [l+1/2]_\q^2 \rightarrow (k+1/2)^2 + (l+1/2)^2.
\end{equation}
These are the eigenvalues   for the commutative spinor Laplacian on the flat square torus with spin structure $\sigma= (1,1) $.
The comparison of the Dirac eigenvalues with the corresponding ones  of the commutative torus \eqref{eq:classevBB} with $\sigma=(1,1)$  is shown in Figure \ref{fig:3dcomp}, where a portion of each surface is shown side-by-side.
 Notice that for small values of $k$ and $l$, the spectra of the fuzzy torus and commutative torus  are practically identical. Physically, this means that one cannot detect the non-commutative behaviour of the fuzzy torus at low energies. 


\begin{figure}
\begin{center}
\includegraphics[scale=0.5]{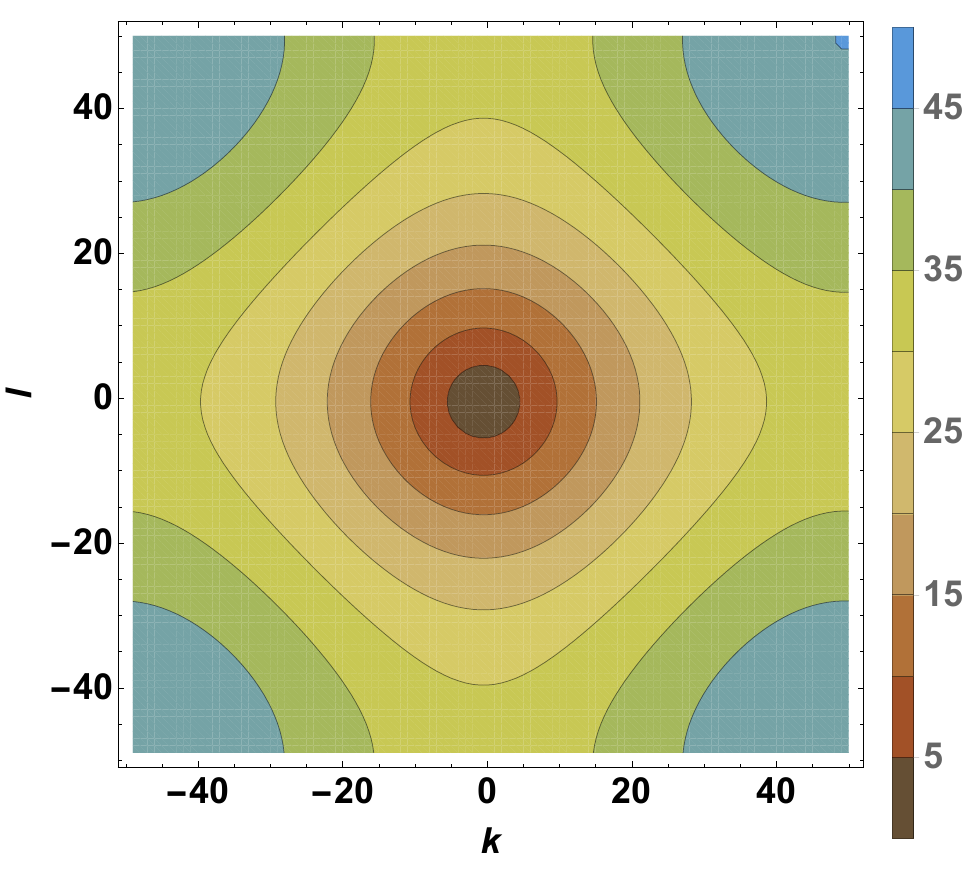}
\end{center}
\vspace{-0.8cm}
\caption{Contour plot for the positive eigenvalues \eqref{eq:fuzzyevBB} with $N=100$.}\label{fig:contoursq}
\end{figure}
It is also useful to represent the spectrum through a contour plot. This will allow an easy comparison of the geometries of different tori. The contour plot for the square fuzzy torus is shown in Figure \ref{fig:contoursq}. The plot has some remarkably symmetrical features that encode information about the geometry and the multiplicities. The  lines on which the large multiplicity eigenvalues lie are defined by $ (l+1/2) \pm (k+1/2) = \pm^{\prime} (N/2 + 2jN) $ with $ j \in \mathbb{Z}$. 

\begin{figure} 
\begin{center}
\includegraphics[scale=0.8]{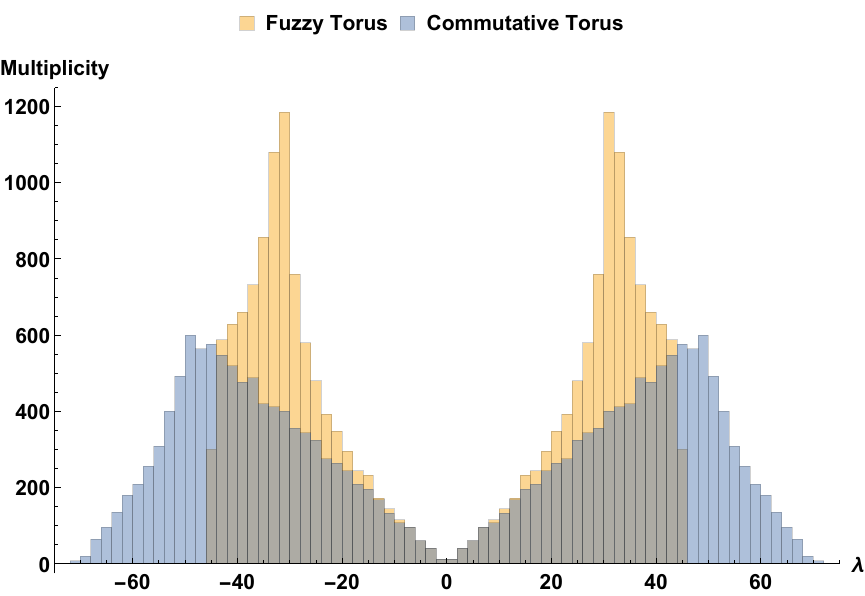}
\end{center}
\caption{Histogram of multiplicity against the eigenvalues \eqref{eq:fuzzyevBB} and  \eqref{eq:commtorusev}, with $N=100$ and bin width $2$.}\label{fig:hist}
\end{figure}

A comparison of the multiplicities of the commutative torus and the fuzzy torus is shown in the histogram in Figure \ref{fig:hist}. The commutative torus multiplicities follow an approximate V shape until the artificial cut-off \eqref{eq:cutoff} is reached. One sees that for small eigenvalues, the fuzzy torus and commutative torus have approximately the same multiplicities but diverge from each other at higher energies.

\subsection{Fuzzy torus with integral metric}\label{sec:fuzzyint}
An integral metric is determined by the choice of $X$ and $Y$.
 Similarly to Section \ref{sec:scalarlaplace}, let  $ \mathfrak h = M_N(\mathbb{C}) $, $J = \ast $, $U=C$ and $V=S$, so that  $\mathcal A= M_N(\C)$. Now set $X=E^{(a,b)},Y = E^{(c,d)}$, so that $ XY = \X YX $ with $ \X^{1/2} = \q^{(ad-bc)/2} $. 
Note that different choices of $U$ and $V$ are possible, as will be discussed in Section \ref{sec:fuzzyspin}. With the choice given here the fuzzy torus is irreducible. 

The action of the group $ \mathbb{Z}_N \times \mathbb{Z}_N$ on the fuzzy torus algebra, as given in Definition \ref{def:translation}, is 
the non-commutative analogue of the group action $U(1) \times U(1) $ on the commutative torus. This action can be extended to the spectral triple. It is shown here that the fuzzy torus Dirac operator is equivariant under this  action and that this is analogous to the explicit formulas for the equivariance of the commutative integral Dirac operator shown in Section \ref{sec:commintegral}.

Let $\q^{1/2}=e^{i\pi K/N}$.  For $j,n\in\Z$, define
\begin{equation}
\Theta = \frac{2 \pi K j}{N}, \quad \Phi = \frac{2 \pi K n}{N}
\end{equation}  
and
\begin{equation}
W = \text{exp}\left(-\frac{1}{2} ( \gamma^1 \gamma^2(  a \Theta + b \Phi) +  \gamma^3 \gamma^4( c \Theta + d \Phi) ) \right).
\end{equation}
The action on the Hilbert space of the spectral triple is 

\begin{equation}
\Pi_{(j,n)} ( v \otimes m ) = W\, v \otimes P_{(j,n)}^{\vphantom{-1}}\, m\, P_{(j,n)}^{-1},
\end{equation}
with $P_{(j,n)}$ given in Definition \ref{def:translation}.
This is the non-commutative analogue of \eqref{eq:pispinaction}.
The fuzzy torus  Dirac operator is equivariant with respect to this action, i.e., 
\begin{equation}
D_{X,Y} \Pi_{(j,n)} =  \Pi_{(j,n)}D_{X,Y}.
\end{equation}

Since the eigenvalues of $\gamma^1\gamma^2$ and $\gamma^3\gamma^4$ are $\pm i$, the action is periodic. If $N$ is odd
\begin{equation}  \Pi_{(N,0)}=1, \quad \Pi_{(0,N)}=1,
\end{equation}
whereas if $N$ is even,
\begin{equation}  \Pi_{(N,0)}=(-1)^{a+c}, \quad \Pi_{(0,N)}=(-1)^{b+d}.
\end{equation}
Thus $\Pi_{(j,n)}$ is a projective representation of $\Z_N\times\Z_N$.

Define the canonical spin structure for this geometry
\begin{equation}\sigma_c=(\lBracket a\rBracket_2+\lBracket c\rBracket_2,\; \lBracket b\rBracket_2+\lBracket d\rBracket_2)= \lBracket A^T\rBracket_2\,(1,1),
\end{equation}
according to \eqref{eq:sspb}.
Then parameters $k$ and $l$ are defined by
$2(k,l) =\sigma_c \pmod 2$, as in \eqref{eq:sskl}. 
These parameters determine vectors
\begin{equation}
\psi = \left( \begin{array}{c}
\alpha\,  \e^{(k+(a+c)/2,\, l+(b+d)/2)}
\\
\beta\,  \e^{(k-(a+c)/2,\,l-(b+d)/2)}
\\
\gamma\, \e^{(k+(c-a)/2,\,l+(d-b)/2)}
\\
\delta\, \e^{(k+(a-c)/2,\,l+(b-d)/2)}
\end{array} \right),\label{eq:eigenvectors}
\end{equation}
with $\alpha,\beta,\gamma,\delta$ arbitrary complex numbers. A short calculation shows these are eigenvectors for the action of the translations, i.e.,
\begin{equation}
\Pi_{(j,n)} \psi = e^{i(\Theta k + \Phi l)} \psi.
\end{equation}

The $\psi$  are also eigenvectors of $D^2_{X,Y}$, the eigenvalue equation being
\begin{equation}
D^2_{X,Y} \psi = \bigg{(} \frac{[al-bk]_\q^2+[dk-cl]_\q^2}{[ad-bc]_\q^2} \bigg{)} \psi
 = \left( \left[ \frac{al-bk}{ad-bc} \right]_\X^2  + \left[ \frac{dk-cl}{ad-bc} \right]_\X^2 \right) \psi.\label{eq:evim}
\end{equation} 
One can find the eigenvectors of $D_{X,Y}$ itself by solving a $4\times 4$ matrix equation for the coefficients $\alpha, \beta, \gamma, \delta$. The corresponding eigenvalues are the two square roots of the eigenvalues in \eqref{eq:evim}, each with multiplicity two. These formulas for the eigenvectors and eigenvalues are direct analogues of the commutative case, for example as shown explicitly for the square torus in Section \ref{sec:squaretorus}.

Now suppose that $k$ and $l$ are both integers, i.e., $\sigma_c=0$. Then the eigenvalues are exactly the same as for the scalar Laplacian in \eqref{eq:laplaceev-integral}. There are in fact \emph{two} natural unitary transformations relating  $1\otimes \Delta_{U,V} $ and $ D^2_{U,V}$. Define algebra elements  $ T = E^{(\frac{a+c}{2},\frac{b+d}{2})} $ and $ H = E^{(\frac{c-a}{2},\frac{d-b}{2})} $. Using the defining (left) action of these operators in $\mathcal H$,
the first unitary transformation is  
\begin{equation}\label{eq:UL}\begin{aligned}
\mathcal{U}_L &= \frac{1}{4} (i(\gamma^1 \gamma^2 - \gamma^3 \gamma^4)(H^*-H)-i(\gamma^1 \gamma^2 + \gamma^3 \gamma^4)(T^*-T) \\ &+(1+\gamma)(H+H^*) + (1-\gamma)(T+T^*)), 
\end{aligned}\end{equation}
while the second one is
\begin{equation}
\mathcal{U}_R =\mathcal J\mathcal{U}_L\mathcal J^{-1}. \end{equation}
Then
\begin{equation}
1\otimes \Delta_{X,Y} = \mathcal{U}_L D^2_{X,Y} \mathcal{U}^*_L = \mathcal{U}_R D^2_{X,Y} \mathcal{U}^*_R.
\end{equation}
 One may think of both $\mathcal{U}_L$ and  $\mathcal{U}_R$  as non-commutative analogues of the spin transformation $Z$ introduced in \eqref{eq:spinlift}. In fact, the formula \eqref{eq:UL} becomes exactly \eqref{eq:spinlift} if $T$ and $H$ are replaced by the corresponding exponential functions of the variables $\theta$ and $\phi$.

 \begin{figure}
\begin{center}
\includegraphics[scale=0.5]{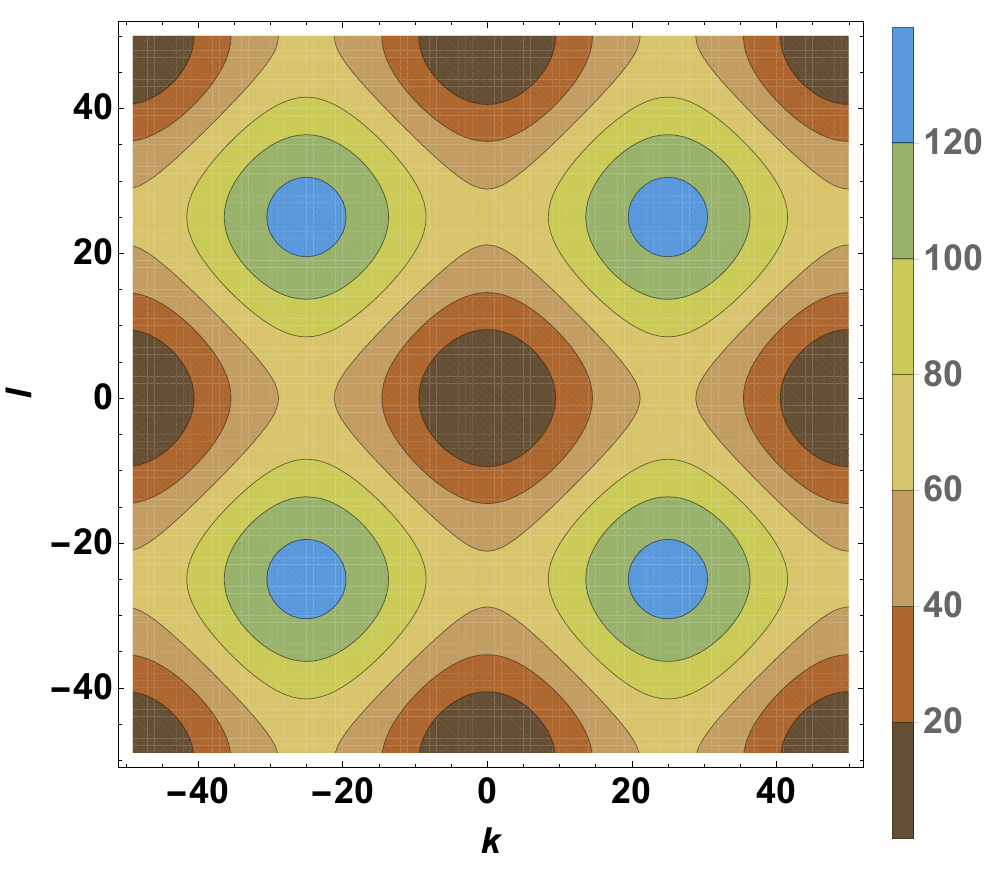}
\end{center}
\vspace{-0.8cm}
\caption{Contour plot for the spectrum \eqref{eq:fuzzyevsq} with $N=100$.}\label{fig:fuzzyevsq}
\end{figure}

\begin{example} Using an instance of the fuzzy torus of Example \ref{ex:genXY} with  $ X = C^2 $ and $ Y = S^2$, the eigenvalues of the Dirac operator are 
\begin{equation}
\lambda_{k,l,\pm} =  \pm \sqrt{\left[\frac{k+1}{2}\right]_\X^2 + \left[\frac{l+1}{2}\right]_\X^2}. \label{eq:fuzzyevsq}
\end{equation}
The eigenvalues are presented in a contour plot in Figure \ref{fig:fuzzyevsq}. It can be seen that the torus has now been stretched equally in both directions by a factor of two. 
\end{example}

\begin{figure}
\begin{center}
\includegraphics[scale=0.5]{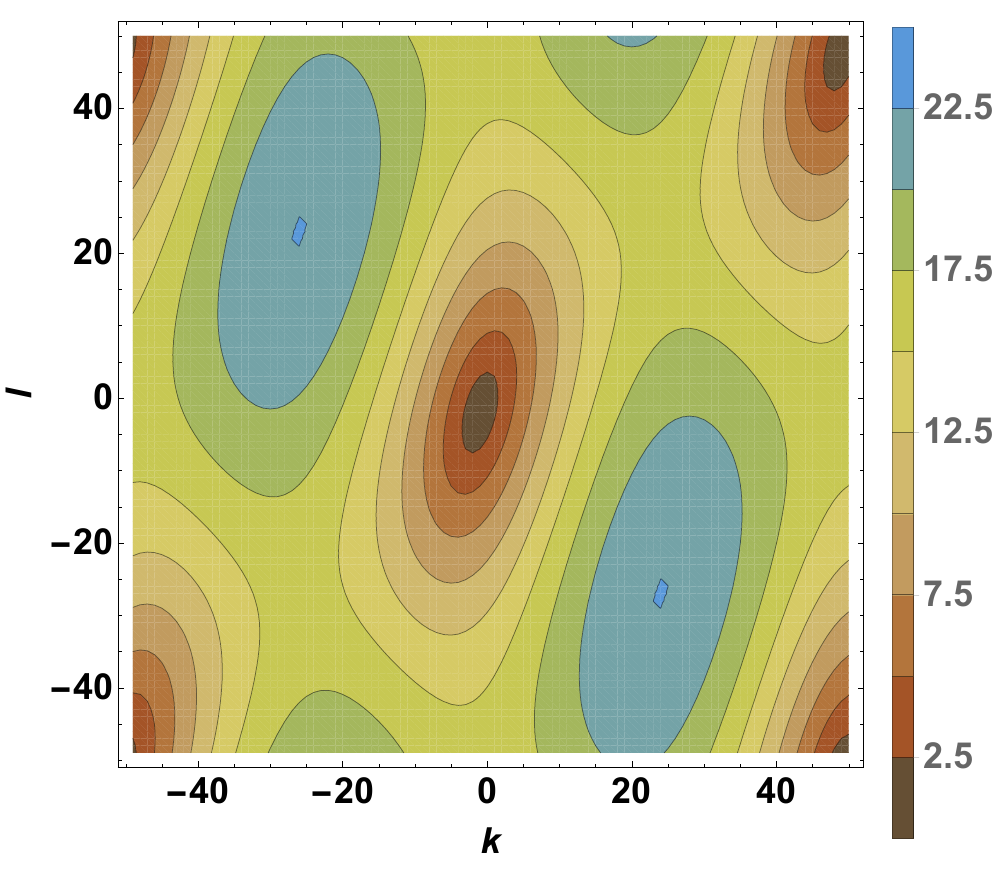}
\end{center}
\vspace{-0.8cm}
\caption{Contour plot for the eigenvalues \eqref{eq:fuzzyevher} and $N=100$.}\label{fig:wonky}
\end{figure}

\begin{example} Another instance of Example \ref{ex:genXY} is with $ X = \q^{-1/2} CS$ and  $Y = S^2 $. The eigenvalues of the Dirac operator are
\begin{equation}
\lambda_{k,l,\pm} =  \sqrt{\left[ \frac{l-k + 1}{2} \right]_\X^2  + \left[ k + \frac12 \right]_\X^2}. \label{eq:fuzzyevher}
\end{equation}
Again, a contour plot of these eigenvalues, in Figure \ref{fig:wonky}, gives a visual representation of the geometry. 
This plot shows clearly the non-square geometry of this torus.
\end{example}

\subsection{Spin Structures on the fuzzy torus}\label{sec:fuzzyspin}

Section \ref{sec:fuzzyint} showed how to construct a non-commutative spectral triple for each integral metric on the torus and a value of $N$. The commutative analogue of this has a particular spin structure, the canonical spin structure.  In this section, it is shown how to construct a spectral triple for a fuzzy torus with any integral metric and any spin structure.  This is done by constructing a non-commutative generalisation of the appropriate covering space of a torus. 

The definitions could be phrased more generally but since it is not yet clear whether there are other interesting examples, this is left as an open problem. There is a very interesting analogy with symplectic reduction, which is outlined at the end of this section, and may help with finding the appropriate non-commutative context.

The commutative construction that is to be generalised is as follows. There is an analogue of the Hurewicz homomorphism for $\Z_2$ coefficients. This is the homomorphism $h\colon\pi_1(\mathcal M)\to H_1(\mathcal M;\Z_2)$ that takes a loop to the corresponding $\Z_2$-cycle. It is obtained from the usual Hurewicz homomorphism $h\colon\pi_1(\mathcal M)\to H_1(\mathcal M;\Z)$ by tensoring with $1\in\Z_2$.
A Riemannian manifold $\mathcal M$  has a regular covering space $\widehat {\mathcal M}$ associated to $h$; this is a principal bundle with group $H_1(\mathcal M;\Z_2)$. Now suppose $\mathcal M$ has spin structure $s$. The spinor bundle $\mathcal W$ on $\mathcal M$ pulls back to a spinor bundle $\widehat {\mathcal W}$ on $\widehat {\mathcal M}$, with a pull-back Clifford multiplication, and hence a uniquely-determined Dirac operator $\widehat D$. Sections of $\widehat {\mathcal W}$ that are equivariant with respect to a character $\chi\colon H_1(\mathcal M;\Z_2)\to\U(1)$, i.e.,
\begin{equation}g\psi=\chi(g)\psi\label{eq:split}\end{equation}
 can be considered sections of a spin bundle on $\mathcal M$ with another spin structure $s'$. This spin structure is determined by $s'=s+\pi^\bullet t$, with $t\in H^1(\mathcal M;\Z_2)$ the cohomology class that corresponds to $\chi$, i.e., $\chi=(-1)^t$. Then $\widehat D$ restricted to these sections gives the Dirac operator on $\mathcal M$ with the spin structure $s'$.
 
\paragraph{Non-commutative coverings}
 The basic structure is an inclusion of $*$-algebras $\mathcal A\subset\mathcal B$. In the commutative case, this is the inclusion of functions that is dual to the projection map of the covering.  

For the non-commutative example, both $\mathcal A$ and $\mathcal B$ are fuzzy tori. It suffices to consider the case where the `total space'  is the fuzzy torus $\bigl(C,S,\mathfrak h=M_{N}(\C)\bigr)$ with clock and shift operators satisfying the relation
 $CS= \q SC$, and with real structure $J=*$. Then $\mathcal B=\langle C,S\rangle=M_{N}(\C)$.

Now suppose $U,V\in\mathcal B$ satisfying $UV=\Q VU$, with $N'$ the order of $\Q$. Then $\mathcal A=\langle U,V\rangle$ is the algebra for the second fuzzy torus $(U,V, \mathfrak h)$ forming the `base space'. The elements $U^{N'}$ and $V^{N'}$ are central in $\mathcal A=\langle U,V\rangle$ and so generate a finite abelian group $\mathcal G$. 

This group acts in $\mathfrak h$ by the adjoint action $\psi\mapsto g\psi g^{-1}$, and as automorphisms of $\mathcal B$ by $b\mapsto gbg^{-1}$. Since elements of the subalgebra $\mathcal A$ are invariant under this action, it is a non-commutative analogue of the deck transformations of a covering. 

The Hilbert space $\mathfrak h$ splits into subspaces according to the unitary characters of $\mathcal G$, i.e., $\mathfrak h_\chi$ is the subspace of vectors $\psi$ satisfying
\begin{equation}g\psi g^{-1}=\chi(g)\psi,\label{eq:ncsplit}\end{equation}
a non-commutative version of \eqref{eq:linesplit}. Then $\mathfrak h=\oplus_\chi \mathfrak h_\chi$ is the non-commutative version of \eqref{eq:hdecomp}.
Acting with $J$ shows
\begin{equation}g(J\psi) g^{-1}=gJg\psi=J(g\psi g^{-1})=J\chi(g)\psi=\overline\chi(g)J\psi,\label{eq:ncsplit2}\end{equation}
so that $J\psi$ lies in $\mathfrak h_{\overline\chi}$, with $\overline\chi$ the complex conjugate character. 

In particular, $J\psi$ lies in $\mathfrak h_\chi$ if $\chi=\overline \chi$, which is true if all elements of $\mathcal G$ have order two.  In this case, each $(U,V,\mathfrak h_\chi)$ is a fuzzy torus.

\paragraph{Non-commutative spin structures} Specialising to the case that is analogous to spin structures, it is necessary to have $U^{N'}\ne 1$, $V^{N'}\ne1$ and the relations $U^{2N'}=V^{2N'}=1$, so that $\mathcal G\cong\Z_2\times\Z_2$. This can be achieved by setting $N=4N'$ and
\begin{equation}U=  C^2, \quad V= S^2
\end{equation}
to define $\mathcal A\subset \mathcal B$. The parameters are related by $\q^4=\Q$. Note that $\mathcal G$ acts non-trivially on the algebra $\mathcal B$: for example, $U^{N'} S U^{-N'}=-S$.

According to Example \ref{ex:C2S2}, $\mathcal A=\mathcal A_{00}\oplus\mathcal A_{01}\oplus\mathcal A_{10}\oplus\mathcal A_{11}$ and, moreover, each summand is a $*$-algebra. There is also a surjective homomorphism $\mathcal A\to\mathcal A_{00}\cong M_{N'}(\C)=\langle C',S'\rangle$ given by $U\mapsto C'$, $V\mapsto S'$.

Spectral triples are constructed according to Definition \ref{def:spectral} using the Dirac operator $D_{X,Y}$, for $X,Y\in\mathcal A$. This defines spectral triples for the fuzzy torus $\bigl(C,S,\mathfrak h\bigr)$, with algebra $\mathcal B$, and also for the fuzzy torus $(U,V, \mathfrak h)$, with algebra $\mathcal A$. Since the characters are all real, $\chi=\overline \chi$, the latter splits into fuzzy tori $(U,V, \mathfrak h_\chi)$, giving spectral triples with Hilbert space $\mathcal H_\chi=\C^4\otimes \mathfrak h_\chi$.

The four characters of $\mathcal G$ are labelled $00, 01, 10, 11$ and the splitting is written 
\begin{equation}\mathcal H=\mathcal H_{00}\oplus\mathcal H_{01}\oplus\mathcal H_{10}\oplus\mathcal H_{11},\end{equation}
 in a similar way to Example \ref{ex:C2S2}, corresponding to the bundles with periodic or anti-periodic boundary conditions in the commutative analogue. Using this analogy, the spin structure of $T^2$ associated with $\mathcal H_\chi$ should be defined as $\sigma=\sigma_c+t$, with $\chi=(-1)^t$.

\begin{example} The simplest example is for the square fuzzy torus determined by $D_{U,V}$.
The eigenvectors of $D_{U,V}$ are given by \eqref{eq:eigenvectors} with $a=d=2$, $b=c=0$. They are indexed by integers $k$ and $l$ and lie in $\mathcal H_{hj}$ if $(k+1,l+1)=(h,j)\mod 2$.

The spectrum of $D_{U,V}$ is readily calculated from \eqref{eq:evim}. The eigenvalue of $D_{U,V}^2$ is
\begin{equation}
 \left[ \frac{l}{2} \right]_\Q^2  + \left[ \frac{k}{2} \right]_\Q^2 .\label{eq:evimsp}
\end{equation} 
One can see that the eigenvalues are exactly the same as in $\eqref{eq:CSev}$ in the case $h=j=0$, so that the spin structure is $\sigma=(1,1)$. However the multiplicity of each eigenvalue in \eqref{eq:evimsp} is four times greater. For general $h, j$, the eigenvalues correspond to the commutative case \eqref{eq:classevBB} if $\sigma=(h+1,j+1)$. Thus one identifies $\sigma_c=(1,1)$ as expected.
\end{example}

\paragraph{Symplectic reduction}
As a final remark, the algebra homomorphisms $\mathcal A_{00}\leftarrow\mathcal A\hookrightarrow \mathcal B$ form a discrete and non-commutative analogue of the maps in a Marsden-Weinstein reduction of a symplectic manifold by a group action \cite{MarsdenWeinstein}. There, starting with a Lie group $G$ and a Hamiltonian $G$-action on a symplectic manifold $\mathcal M$ with momentum map $\mu$, one has the maps $\mathcal M\git G\hookrightarrow \mathcal M/ G\leftarrow \mathcal M$, with  $\mathcal M\git  G=\mu^{-1}(0)/ G$.

In the non-commutative example, the analogue of $G$ is the finite group  $\mathcal G$.
The algebra $\mathcal B$ is the non-commutative analogue of (functions on) $\mathcal M$, $\mathcal A$ is the analogue of the space of orbits $\mathcal M/ G$, while $\mathcal A_{00}$ is the analogue of the symplectic quotient $\mathcal M\git G$ .

 Denote the group algebra of $\mathcal G$ by $\C[\mathcal G]$. This is dual to the group of unitary characters $\mathcal G^*$, i.e., $\C[\mathcal G]\cong C(\mathcal G^*)$. The inclusion $\C[\mathcal G]\hookrightarrow \mathcal B$ is the analogue of the momentum map $\mu$, the (discrete) momenta being the elements of $\mathcal G^*$. Indeed, the
 `constraint space' $\mathcal A_{00}$  corresponds exactly to momentum $1\in\mathcal G^*$.

\end{document}